\documentclass[final]{siamltex}
\usepackage{graphicx, amssymb, amsmath, epsfig}
\usepackage{color}
\usepackage{showlabels}
\usepackage{hyperref}
\usepackage{booktabs}
\usepackage{boxedminipage}
\usepackage[ruled]{algorithm2e}

\graphicspath{{./}{figs/}}
\newcommand{\ra}[1]{\renewcommand{\arraystretch}{#1}}

\newcommand{\R}{\mathbb R}

\newcommand{\A}{\,\mathcal A\,}
\newcommand{\T}{\,\mathcal T\,}

\begin{document}

\title{An Efficient  Policy Iteration Algorithm for Dynamic Programming Equations \thanks{The authors wish to acknowledge the support obtained by the following grants: AFOSR Grant no.\ FA9550-10-1-0029,
ITN - Marie Curie Grant n. 264735-SADCO , START-Project Y432-N15 {``Sparse Approximation and
Optimization in High-Dimensions''} and SAPIENZA 2009 {``Analisi ed approssimazione di modelli differenziali nonlineari in fluidodinamica e scienza dei materiali''}. The authors also wish to thank the CASPUR Consortium for its technical support.}
}

\author{
Alessandro Alla\footnote{Department of Mathematics,
Universit\"at Hamburg, Bundesstr. 55, 20146 Hamburg, Germany.
\texttt{alessandro.alla@uni-hamburg.de}},
Maurizio Falcone\footnote{(corresponding author) Dipartimento di Matematica,
SAPIENZA - Universit\`a di Roma, P.le Aldo Moro 2, 00185 Rome, Italy.
\texttt{falcone@mat.uniroma1.it}},
Dante Kalise\footnote{Johann Radon Institute for Computational and Applied Mathematics,
Altenbergerstra\ss e 69, 4040 Linz, Austria.
\texttt{dante.kalise@oeaw.ac.at}}
}

\maketitle

\begin{abstract}\noindent
We present an accelerated algorithm for the solution of static Hamilton-Jacobi-Bellman equations related to optimal control problems. Our scheme is based on a classic policy iteration procedure, which is known to have superlinear convergence in many relevant cases provided the initial guess is sufficiently close to the solution. This limitation often degenerates into a behavior similar to a value iteration method, with an increased computation time. The new scheme circumvents this problem by combining the advantages of both algorithms with an efficient coupling. The method starts with a coarse-mesh value iteration phase and then switches to a fine-mesh policy iteration procedure when a certain error threshold is reached. A delicate point is to determine this threshold in order to avoid cumbersome computations with the value iteration and, at the same time, to be reasonably sure that the policy iteration method will finally converge to the optimal solution.
We analyze the methods and efficient coupling in a number of examples in dimension two, three and four illustrating their properties.
\end{abstract}



\pagestyle{myheadings}
\thispagestyle{plain}


\section{Introduction}
The numerical solution of optimal control problems is a crucial issue for many industrial applications such as aerospace engineering, chemical processing, power systems, and resource economics, among many others. In some cases the original problem comes from a different setting, e.g. when one has to fit a given set of data or has to solve a shape optimization problem, but has been reformulated in terms of a control problems for an appropriate dynamic and cost functional. The typical goal is then to compute an optimal trajectory for the controlled system and the optimal control corresponding to it. In the framework of open-loop controls the classical solution is based on the Pontryagin Maximum Principle which leads to the solution of a two-point boundary value problem for the coupled state/costate system. The numerical solution can be obtained via a shooting method. Despite its simplicity and mathematical elegance, this approach is not always satisfactory because the initialization of the shooting method can be a difficult task, mainly for the costate variables. Moreover, this approach is typically based on necessary conditions for optimality and produces only open-loop controls.
It is well known that the Dynamic Programming (DP) approach introduced by Bellman \cite{B57} produces optimal control in feedback form so it looks more appealing in terms of online implementations and robustness. However, the synthesis of feedback controls require the previous knowledge of the value function and this is the major bottleneck for the application of DP. In fact  the value of optimal control problems are known to be only Lipschitz continuous even when the data are regular and the characterization of the value function is obtained in terms of a first order nonlinear Hamilton-Jacobi-Bellman (HJB) partial differential equation. In the last twenty years,  the DP approach has been pursued for all the classical control problems in the framework of viscosity solution introduced by Crandall and Lions in the 80's (see \cite{BCD97} for a comprehensive illustration of this approach). Moreover, several approximation schemes have been proposed for this class of equations, ranging from finite differences to semi-Lagrangian and finite volume methods. Some of these algorithms converge to the value function but their convergence is slow. The so-called \textsl{curse of the dimensionality}, namely, the fact that the dimension of the partial differential equation characterizing the value function increases as the dimension of the state space does, constitutes a major computational challenge towards a practical implementation of numerical algorithms for optimal control design based on viscosity solutions of HJB equations.

\noindent Our main contribution in this paper is a new accelerated algorithm which can produce an accurate approximation of the value function in a reduced amount of time in comparison to the currently available methods. Furthermore, the proposed scheme can be used in a wide variety of problems connected to static HJB equations, such as infinite horizon optimal control, minimum time control and some cases of pursuit-evasion games. The new method couples two ideas already existing in the literature: the value iteration method (VI) and the policy iteration method (PI) for the solution of Bellman equations. The first is known to be slow but convergent for any initial guess, while the second is known to be fast when it converges (but if not initialized correctly, convergence might be as slow as the value iteration). The approach that we consider relates to multigrid methods (we refer to Santos \cite{SA99} for a brief introduction to subject in this context), as the coupling that we introduce features an unidirectional, two-level mesh. The work by Chow  and Tsitsiklis \cite{CT91} exploits a similar idea with a value iteration algorithm. However, as far as we know the efficient coupling between the two methods has not been investigated.

\noindent To set this paper into perspective, we must recall that algorithms based on the iteration in the space of controls (or policies) for the solution of HJB equations has a rather long history, starting more or less at the same time of dynamic programming. The PI method, also known as Howard's algorithm \cite{H60},  has been investigated by Kalaba \cite{K59} and  Pollatschek and Avi-Itzhak \cite{PA}    who proved  that  it corresponds to the Newton method applied to the functional equation of dynamic programming. Later, Puterman and Brumelle \cite{PB79} have given sufficient conditions for the rate of convergence to be either superlinear or quadratic. More recent contributions on the policy iteration method and some extensions to games can be found in Santos and Rust \cite{SR04} and  Bokanowski et al. \cite{BMZ09}. Results on its numerical implementation and diverse hybrid algorithms related to the proposed scheme have been reported in Capuzzo-Dolcetta and Falcone \cite{CDF89}, Gonz\'alez and Sagastiz\'abal \cite{GS90}, Gr\"une \cite{G97} and in the recent monograph by Falcone and Ferretti \cite{FF14}.

\noindent Finally, we should mention that an acceleration method based on the the set of subsolutions has been studied in Falcone \cite{F87} (see also Tidball and Gonz\'alez \cite{TG92} for a specific application to the Isaacs equation).
More in general, dealing with acceleration methods for Hamilton-Jacobi-Bellman equations, we should also mention approaches based on domain decomposition algorithms as in Falcone et al. \cite{FLS94} and more recently by Cacace et al. \cite{CCFP12}, on geometric considerations as in Botkin, et al. \cite{BHT11}, and those focusing on the localization of numerical schemes which leads to Fast Marching Methods. This approach has shown to be very effective for level-set equations related to front propagation problems (see e.g. the book by Sethian \cite{S99}), i.e. eikonal type equations.   At every iteration, the scheme is applied only on a subset of nodes (localization) which are the nodes close to the front, the so-called \textsl{narrow band}.  The remaining part of the grid is divided into two parts: the accepted region, where the solution has been already computed, and the far region where the solution will be computed little by little in the following iterations. At every iteration, one node is accepted and moved from the narrow band to the accepted region; the narrow band is then updated adding the first neighbors of that node (which before where in the far region). For eikonal type equations these methods converge in finite number of iterations to the correct viscosity solution and have a very low complexity (typically $O(N \ln(N)$) where N is the cardinality of the grid). More recently several efforts have been made to extend these methods to more complex problems where the front propagation is anisotropic \cite{SV03} and/or to more general Hamilton-Jacobi equations as in  \cite{AM12}.  However, their implementation is rather delicate and their convergence to the correct viscosity solution  for general Hamilton-Jacobi equations is still an open problem; we refer to \cite{CCF13} for a an extensive discussion and several examples of these limitations.

\noindent The paper is organized as follows. In Section 2, we introduce some basic notions for optimal control synthesis by the dynamic programming principle. Section 3 contains the core of the proposed accelerated method and discuss practical implementation details. In Section 4 we develop a convergence estimate for the proposed scheme. Finally, Section 5 shows our numerical results on  a number of different examples concerning infinite horizon optimal control, minimum time control, and some further extensions towards the optimal control of partial differential equations. In these series of tests we discuss several properties of the proposed scheme and perform comparisons with the different techniques presented in the article.

\section{Dynamic programming in optimal control and the basic solution algorithms}\label{sec2}
In this section we will summarize the basic results for the two methods as they  will constitute the building blocks for our new algorithm. The essential features will be briefly sketched, and more details can be found in the original papers and in some monographs, e.g. in the classical books by Bellman \cite{B57}, Howard \cite{H60} and for a  more recent setting in framework of viscosity solutions in  \cite{CDF89} and  \cite{BCD97}.\\
Let us first present the method for the classical
{\it infinite horizon problem}. Let the dynamics be given by
\begin{eqnarray}\label{eq:dynconaut}
	\left\{\begin{array}{l}
	\dot{y}(t)=  f(y(t),\alpha(t))\cr
	y(0)=  x
	 \end{array}\right.
\end{eqnarray}
where  $y\in \R^d$, $\alpha \in \R^m$ and $\alpha\in\A\equiv\{ a:\R_+\rightarrow A,\,\text{measurable}\}$. If $f$ is Lipschitz continuous with respect to the state variable and continuous with respect to  $(y,\alpha)$, the classical assumptions for the existence and uniqueness result for the Cauchy problem \eqref{eq:dynconaut} are satisfied. To be more precise, the Carath\'eodory theorem (see \cite{FR75} or  \cite{BCD97}) implies that for any  given control $\alpha(\cdot)\in \A$ there exists a unique trajectory $y(\cdot; \alpha)$ satisfying  \eqref{eq:dynconaut}  almost everywhere. Changing the control policy the trajectory will change and we will have a family of  infinitely many solutions of the controlled system \eqref{eq:dynconaut} parametrized with respect to $\alpha$.
\\
Let us introduce the {\em cost functional} $J:\A\rightarrow\R $ which will be used to select the ``optimal trajectory''. For infinite horizon problem the functional is
\begin{eqnarray}\label{eq:joi}
	J_x(\alpha(\cdot))=\int_0^{\infty} g(y(s),\alpha(s))e^{-\lambda s}ds\,,
\end{eqnarray}
where $g$ is Lipschitz continuous in both arguments and $\lambda>0$ is a given parameter. The function $g$ represents the running cost and $\lambda\in\R_+$ is the discount factor which  allows to compare the costs at different times rescaling the costs at time 0. From the technical point of view, the presence of the discount factor guarantees that the integral is finite whenever $g$ is bounded, i.e.  $||g||_{\infty}\leq M_g$. Let us define the value function of the problem as
\begin{equation}
v(x)=\inf_{\alpha(\cdot)\in\A }J_x(\alpha(\cdot))\,.
\end{equation}
It is well known that passing to the limit in the Dynamic Programming Principle one can obtain a characterization of the value function in terms of the following first order non linear Bellman equation
\begin{equation}\label{eq:hjb}
\lambda v(x)+\max_{a\in A} \{-f(x,a)\cdot Dv(x)- g(x,a)\}=0, \quad \hbox{ for }x\in\R^d\,.
\end{equation}
Several approximation schemes on a fixed grid $G$ have been proposed for \eqref{eq:hjb}. Here we will use a semi-Lagrangian approximation based on a Discrete Time Dynamic Programming Principle. This leads to
\begin{equation}\label{hjbh}
v_{\Delta t}(x) =\min_{a\in A} \{ e^{-\lambda \Delta t} v_{\Delta t}\left(x+\Delta t f\left(x,a\right)\right)+ \Delta t g\left(x,a\right)\}\,,
\end{equation}
where $v_{\Delta t}(x)$ converges to $v(x)$ when $\Delta t\rightarrow 0$.
A natural way to solve \eqref{hjbh} is to write it in fixed point form
\begin{equation}\label{hjbhk}
V_i=\min_{a\in A} \{ e^{-\lambda \Delta t} I[V]\left(x_i+\Delta t f\left(x_i,a\right)\right)+ \Delta t g\left(x_i,a\right)\}\,,\quad i=1,\ldots,N_G
\end{equation}
where $\{x_i\}_{i=1}^{N_G}$ are the grid nodes, $V_i$ is the approximate value for $v(x_i)$  and $I[V]:\R^d\rightarrow \R$ represents an interpolation operator defining, for every point $x$, the polynomial reconstruction based on the values $V_i$
 (see \cite[Appendix A]{BCD97} for more details). Finally, one obtains the following algorithm:

\vskip 1mm
\begin{algorithm}[H]\label{viter}\caption{Value Iteration for infinite horizon optimal control \textbf{(VI)}}
 \SetAlgoLined
 \vskip 1mm
 \KwData{Mesh $G$, $\Delta t$, initial guess $V^0$, tolerance $\epsilon$.}
 \vskip 1mm
 \While{$||V^{k+1}-V^k||\geq \epsilon$}{
  \ForAll{$x_i\in G$}
  {\vskip -4mm\begin{equation}\label{eq:viter}
   \hbox{Solve: } \;V_i^{k+1} =\min_{a\in A} \{ e^{-\lambda \Delta t} I\left[V^k\right]\left(x_i+\Delta t f\left(x_i,a\right)\right)+ \Delta t  g\left(x_i,a\right)\}
  \end{equation}\vskip -4mm}
  $k=k+1$}
\end{algorithm}
\vskip 1mm

\noindent Here $V^k_i$ represents the values at a node $x_i$ of the grid at the $k$-th iteration and $I$ is an interpolation operator acting on the values of the grid; without loss of generality, throughout this paper we will assume that the numerical grid $G$ is a regular equidistant array of points with mesh spacing denoted by $\Delta x$, and we consider a multilinear interpolation operator. Extensions to nonuniform grids and high-order interpolants can be performed in a straightforward manner.

\noindent Algorithm \ref{viter} is referred in the literature as the {\em value iteration method}  because, starting from an initial guess $V^0$,  it modifies the values on the grid according to the nonlinear rule \eqref{eq:viter}. It is well-known that the convergence of the value iteration  can be very slow, since the contraction constant $e^{-\lambda \Delta t}$ is close to 1  when $\Delta t$ is close to 0. This means that a higher accuracy will also require more iterations.  Then, there is a need for an acceleration technique in order to cut the link between accuracy and complexity of the value iteration.

\noindent For sake of clarity, the above framework has  been presented for the infinite horizon optimal control problem. However, similar ideas can be extended to other classical control problems with small changes. Let us mention how to deal with the minimum time problem which we will use in the final section on numerical tests.\\
In the minimum time problem one has to drive the controlled dynamical system \eqref{eq:dynconaut} from its initial state to a given target $\mathcal{T}$. Let us assume that the target is a compact subset of  $\R^d$ with non empty interior and piecewise smooth boundary.
The major difficulty dealing with this problem is that the time of arrival to the target starting from the point $x$
\begin{equation}
t(x,\alpha(\cdot)):=
\left\{
\begin{array} {ll}
\inf\limits_{\alpha\in\A} \{t\in\R_+: y(t,\alpha(\cdot))\in \T\} & \hbox{if } y(t,\alpha(t))\in \T \hbox{ for some }t,\,  \\
+\infty & \hbox{ otherwise, }\\
\end{array}
\right.
\end{equation}
can be infinite at some points. As a consequence, the minimum time function defined as
\begin{equation}
T(x)=\inf\limits_{\alpha\in\A}  t(x,\alpha(\cdot))
\end{equation}
is not defined everywhere unless some controllability assumptions are not introduced. In general, this is a free boundary problem where one has to determine at the same time, the couple $(T, \Omega)$, i.e. the minimum time function and its domain. Nevertheless, by applying the Dynamic Programming Principle and the so-called Kruzkhov transform

\begin{equation}\label{eq:Kruz}
v(x)\equiv
\left\{\begin{array}{ll}
1-\exp(-T(x))&  \hbox{ for } T(x)<+\infty\\
1& \hbox{ for }T(x)=+\infty
\end{array}\right.
\end{equation}

\noindent the minimum time problem is characterized in terms of the unique viscosity solution of the BVP
\begin{equation}\label{eq:Hjb1}
\left\{\begin{array}{cc}
v(x)+\sup\limits_{a\in\mathcal{A}}\{-f(x,a)\cdot Dv(x)\}=1& \hbox{ in }\mathcal{R}\backslash\mathcal{T}\\
v(x)=0 & \hbox{on }\partial \mathcal{T}\,,
\end{array}\right.
\end{equation}
where $\mathcal{R}$ stands for the set of point in the state space where the time of arrival is finite. Then, the application of the semi-Lagrangian method presented for the infinite horizon optimal control problem together with a value iteration procedure leads to following iterative scheme:

\vskip 1mm
\begin{algorithm}[H]\label{vitermtp}\caption{Value Iteration for minimum time optimal control \textbf{(VI)}}
 \SetAlgoLined
 \vskip 1mm
 \KwData{Mesh $G$, $\Delta t$, initial guess $V^0$, tolerance $\epsilon$.}
 {{\bf Set:} $V_i=0$, for all $x_i\in G\cap\T$}
 \vskip 1mm
 \While{$||V^{k+1}-V^k||\geq \epsilon$}{
  \ForAll{$x_i\in G\setminus \T$}
  {\vskip -4mm\begin{equation}
\hbox{Solve: }V_i^{k+1} =\min_{a\in A} \{ e^{-\Delta t} I\left[V^k\right]\left(x_i+\Delta t f\left(x_i,a\right)\right)+ 1-e^{-\Delta t}\}\,
\end{equation}\vskip -4mm}
  $k=k+1$}
\end{algorithm}
\vskip 1mm
\noindent The numerical implementation is completed with the boundary conditions $v(x)=0$ at $\partial\mathcal{T}$ (and inside the target as well), and with $v(x)=1$ at other points outside the computational domain (we refer the reader to \cite{BF90a} for more details on the approximation of minimum time problems).

\noindent\textbf{Policy iteration.} We now turn our attention to an alternative solution method for discretized HJB equations of the form \ref{hjbhk}. The {\em approximation in the policy space} (or policy iteration), and is based on a linearization of the Bellman equation. First, an initial guess for the control for every point in the state space is chosen. Once the control has been fixed, the Bellman equation becomes linear (no search for the minimum in the control space is performed), and it is solved as an advection equation. Then, an updated policy is computed and a new iteration starts. Let us sketch the procedure for the scheme related to the infinite horizon problem.

\vskip 1mm
\begin{algorithm}[H]\label{piter}\caption{Policy Iteration for infinite horizon optimal control \textbf{(PI)}}
 \SetAlgoLined
 \vskip 1mm
 \KwData{Mesh $G$, $\Delta t$, initial guess $V^0$, tolerance $\epsilon$.}
 \While{$||V^{k+1}-V^k||\geq \epsilon$}{\vskip 1mm
  \textsl{Policy evaluation step}:\vskip 1mm
  \ForAll{ $x_i\in G$}
  {\vskip -4mm\begin{equation}\label{policy_ev}
V_i^k = \Delta t g\left(x_i,a_i^k\right)+e^{-\lambda \Delta t} I\left[V^k\right]\left(x_i+\Delta t f\left(x_i,a_i^k\right)\right)
\end{equation}\vskip -4mm}

  \textsl{ Policy improvement step}:\vskip 1mm

  \ForAll{ $x_i\in G$}
  {\vskip -4mm\begin{equation}\label{eq:polup}
a_i^{k+1} =  \arg\min_a \left\{ \Delta t g(x_i,a)+e^{-\lambda \Delta t} I\left[V^k\right](x_i+\Delta t f(x_i,a))\right\}
\end{equation}\vskip -4mm}

  $k=k+1$}
  \vskip 1mm

\end{algorithm}

\vspace{0.3cm}

\noindent Note that the solution of \eqref{policy_ev} can be obtained either by a linear system (assuming a linear interpolation operator) or as the limit
\begin{equation}
V^k=\lim_{m\rightarrow +\infty} V^{k,m}\,,
\end{equation}
of the linear time-marching scheme
\begin{equation}\label{marching}
V_i^{k,m+1} = \Delta t g\left(x_i,a_i^k\right)+e^{-\lambda \Delta t} I\left[V^{k,m}\right]\left(x_i+\Delta t f\left(x_i,a_i^k\right)\right).
\end{equation}
Although this scheme is still iterative, the lack of a minimization phase makes it faster than the original value iteration. \\
The sequence $\{V^k\}$ turns out to be monotone decreasing at every node of the grid. In fact, by construction,
 \begin{eqnarray*}
 V_i^k & = & \Delta t g\left(x_i,a_i^k\right)+e^{-\lambda \Delta t} I\left[V^k\right]\left(x_i+\Delta t f\left(x_i,a_i^k\right)\right) \ge \\
 & \ge & \min_a \left\{\Delta t g(x_i,a)+e^{-\lambda \Delta t} I\left[V^k\right](x_i+\Delta t f(x_i,a))\right\} = \\
 & = & \Delta t g\left(x_i,a_i^{k+1}\right)+e^{-\lambda \Delta t} I\left[V^k\right]\left(x_i+\Delta t f\left(x_i,a_i^{k+1}\right)\right) = \\
 & = & V_i^{k+1}
 \end{eqnarray*}
At a theoretical level, policy iteration can be shown to be equivalent to a Newton method, and therefore, under appropriate assumptions, it converges with quadratic speed. On the other hand, convergence is local and this may represent a drawback with respect to value iterations.

 \vspace{0.2cm}\noindent

\section{An accelerated policy iteration algorithm with smart initialization}
\begin{figure}[ht]
\centering
\begin{tabular}{cc}
\epsfig{file=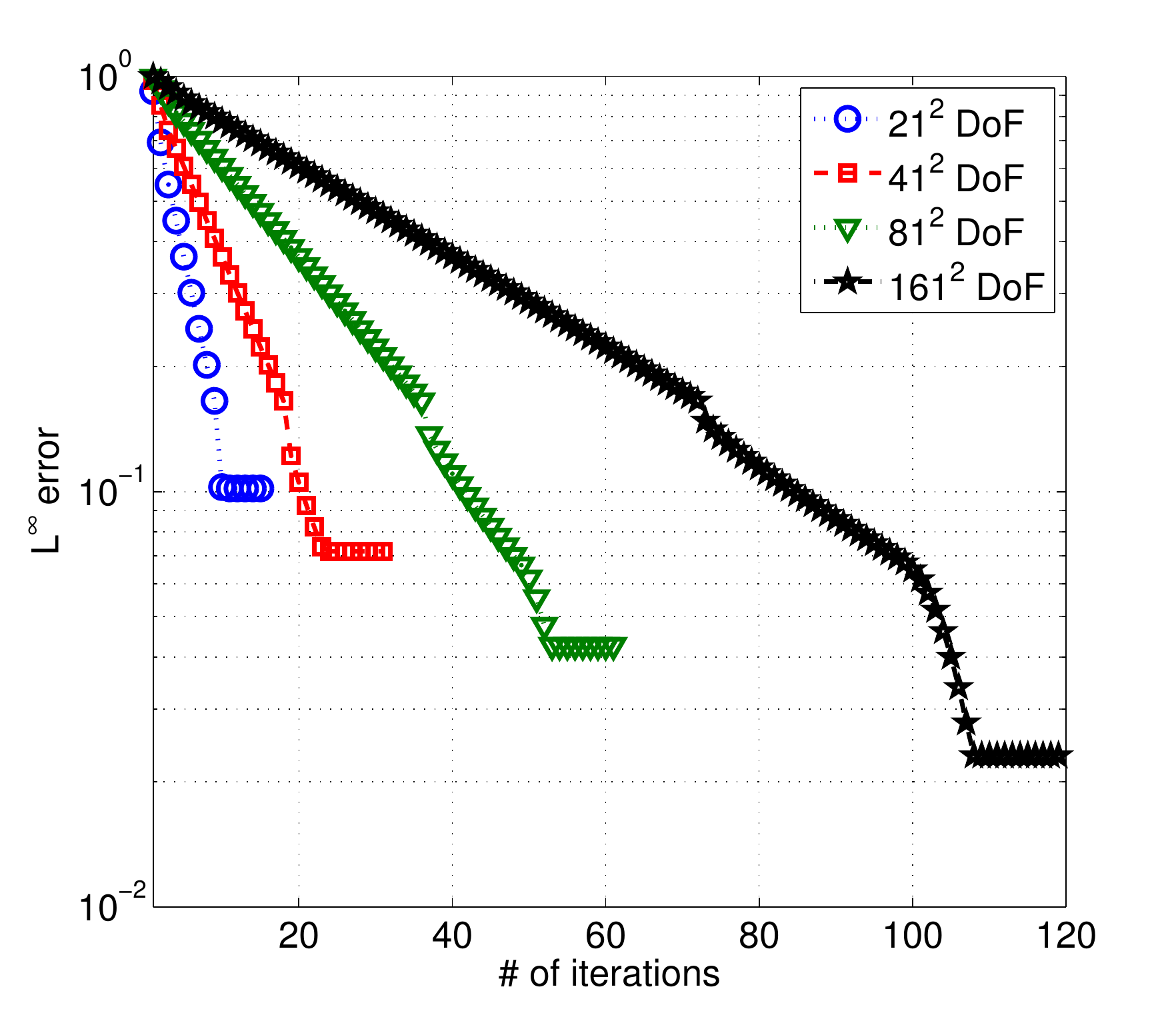,scale=0.35} &
\epsfig{file=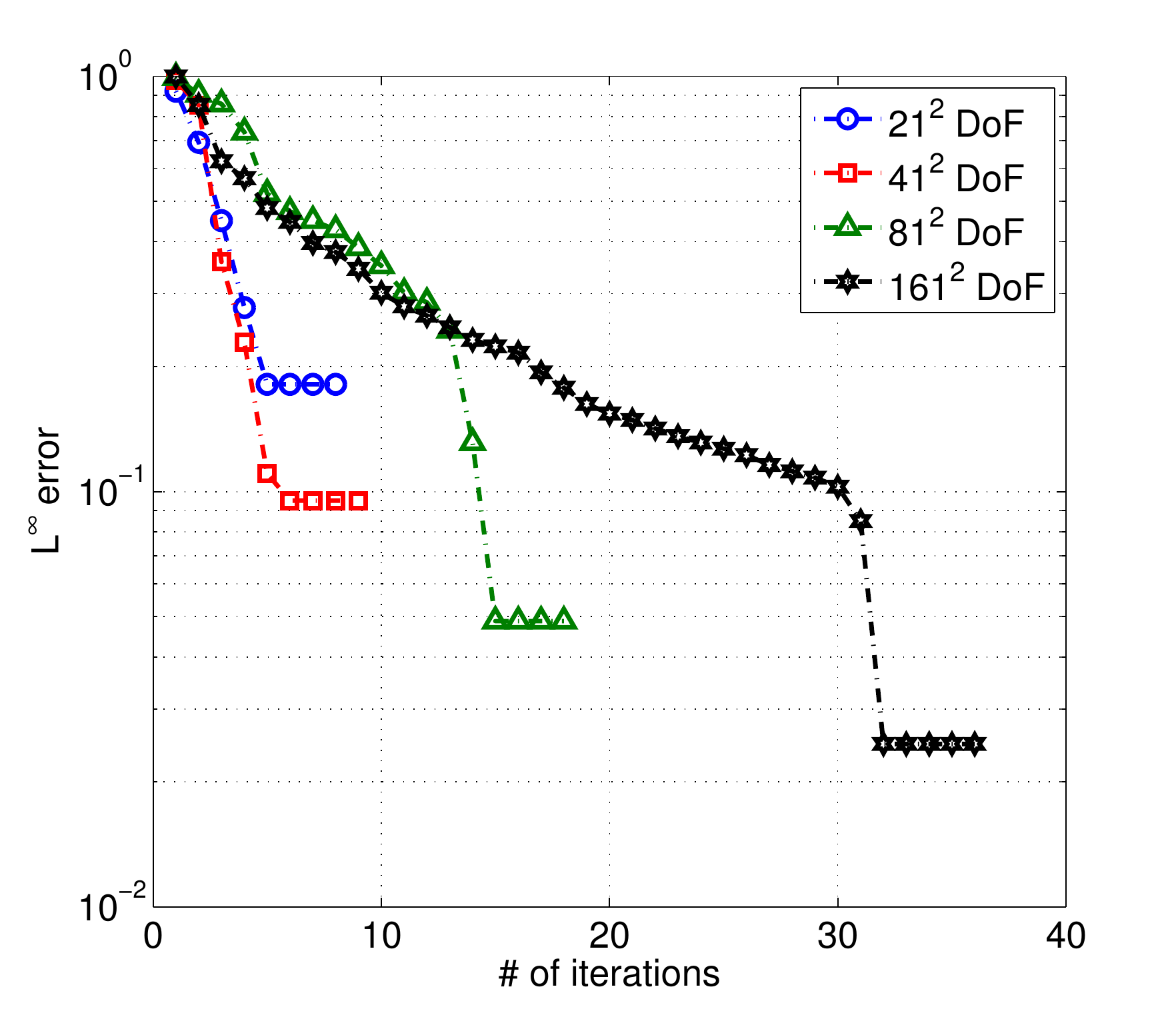,scale=0.35}
\end{tabular}
\caption{Error evolution in a 2D problem: value iteration (left) and policy iteration (right).}
\label{fig:comparison3}
\end{figure}

\begin{figure}[ht]
\centering
\begin{tabular}{cc}
\epsfig{file=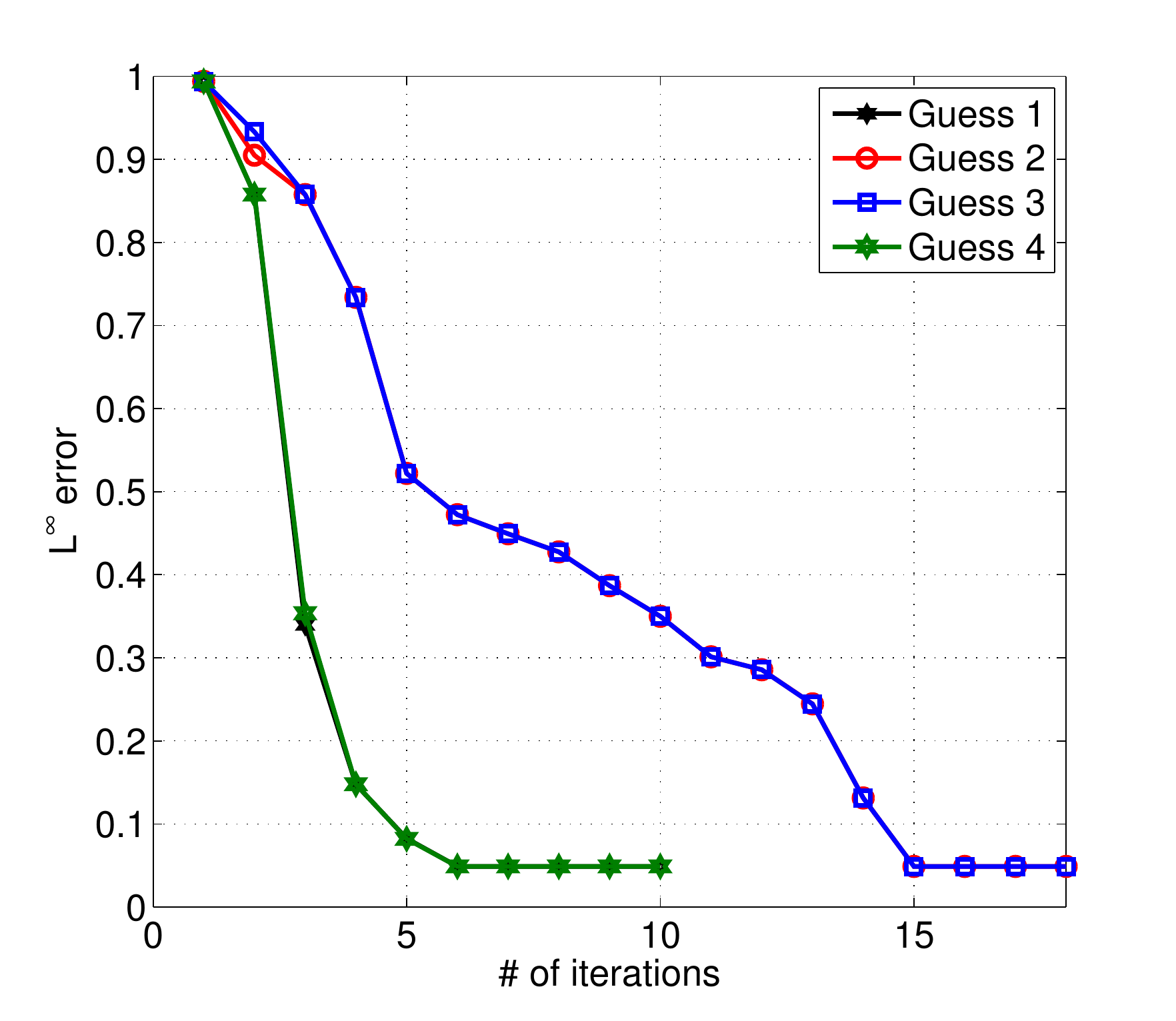,scale=0.35,clip=} &
\epsfig{file=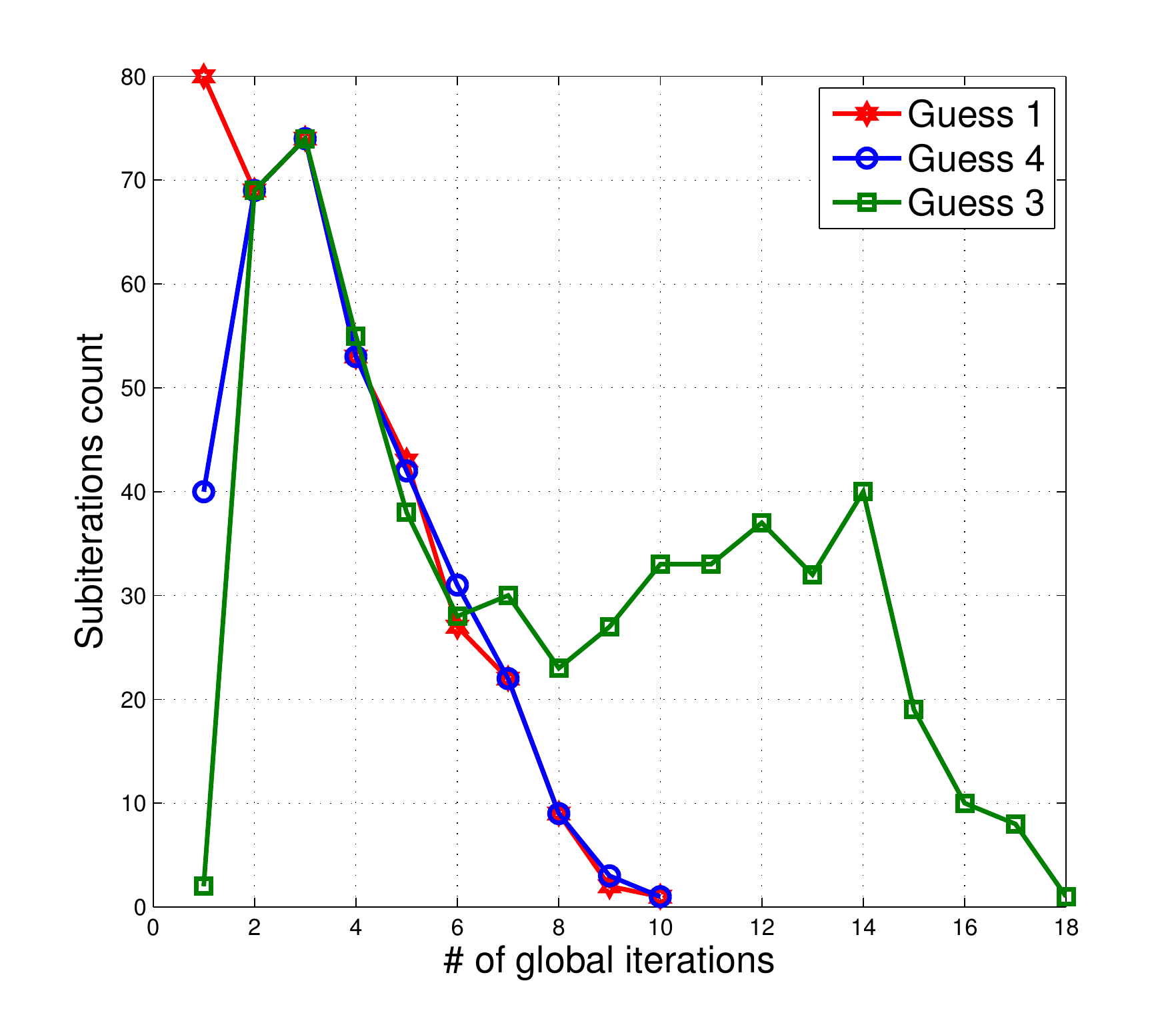,scale=0.35,clip=}
\end{tabular}
\caption{Left: error evolution in a PI algorithm for different initial guesses. Right: evolution of the (sub)iteration count in \eqref{marching}  for different guesses.}
\label{fig:comparison3_2}
\end{figure}

In this section we present an accelerated iterative algorithm which is constructed upon the building blocks previously introduced. We aim at an efficient formulation exploiting the main computational features of both value and policy iteration algorithms. As it has been stated in \cite{PB79}, there exists a theoretical equivalence between both algorithms, which guarantees a rather wide convergence framework. However, from a computational perspective, there are significant differences between both implementations. A first key factor can be observed in Figure \ref{fig:comparison3}, which shows, for a two-dimensional minimum time problem (more details on the test can be found in section 4.4), the typical situation arising with the evolution of the error measured with respect to the optimal solution, when comparing value and policy iteration algorithms. To achieve a  similar error level, policy iteration requires considerable fewer iterations than the value iteration scheme, as quadratic convergent behavior is reached faster for any number of nodes in the state-space grid.
Despite the observed computational evidence, a second issue is observed when examining the policy iteration algorithm in more detail. That is, as shown in Figure \ref{fig:comparison3_2}, the sensitivity of the method with respect to the choice of the initial guess of the control field. It can be seen that different initial admissible control fields can lead to radically different convergent behaviors. While some guesses will produce quadratic convergence from the beginning of the iterative procedure, others can lead to an underperformant value iteration-like evolution of the error. This latter is computationally costly, because it translates into a non-monotone evolution of the subiteration count of the solution of equation \eqref{policy_ev} (if an iterative scheme as in \eqref{marching}).

\noindent A final relevant remark goes back to Figure \ref{fig:comparison3}, where it can be observed that for coarse meshes, the value iteration algorithm generates a fast error decay up to a higher global error. These, combined with the fact that value iteration algorithms are rather insensitive to the choice of the initial guess for the value function (see \cite{SA99} for a detailed error quantification), are crucial points for the construction of our accelerated algorithm. The accelerated policy iteration algorithm is based on a robust initialization of the policy iteration procedure via a coarse value iteration which will yield to a good guess of the initial control field.

\vskip 1mm
\begin{algorithm}[H]\caption{Accelerated Policy Iteration \textbf{(API)}}
 \SetAlgoLined
 \KwData{Coarse mesh $G_c$ and $\Delta t_c$ , fine mesh $G_f$ and $\Delta t_f$, initial coarse guess $V_c^0$, coarse-mesh tolerance $\epsilon_c$, fine-mesh tolerance $\epsilon_f$.}
 \Begin{Coarse-mesh value iteration step: perform Algorithm \ref{viter}\\
 \KwIn{$G_c$, $\Delta t_c$, $V_c^0$, $\epsilon_c$}
 \KwOut{$V_c^*$}
 \ForAll{$x_i\in G_f$}
 {$V^0_f(x_i)=I_1[V^*_c](x_i)$
  $A^0_f(x_i)=\underset{a\in A}{argmin}\;\{e^{-\lambda \Delta t} I_1[V^0_f](x_i+f(x_i,a))+ \Delta t  g\left(x_i,a\right)\}$
 }
Fine-mesh policy iteration step: perform Algorithm \ref{piter}\\
\KwIn{$G_f$, $\Delta t_f$, $V^0_f$, $A^0_f$, $\epsilon_f$}
\KwOut{$V^*_f$}
}
\end{algorithm}
\vskip 1mm
\subsection{Practical details concerning the computational implementation of the algorithm}

The above presented accelerated algorithm can lead to a considerably improved performance when compared to value iteration and naively initialized policy iteration algorithms. However, it naturally contains trade-offs that need to be carefully handled in order to obtain a correct behavior. The extensive numerical tests performed in section 4 suggest the following guidelines:\\

\begin{figure}[ht]
\centering
\scalebox{.99}{
\begin{tabular}{cc}
\epsfig{file=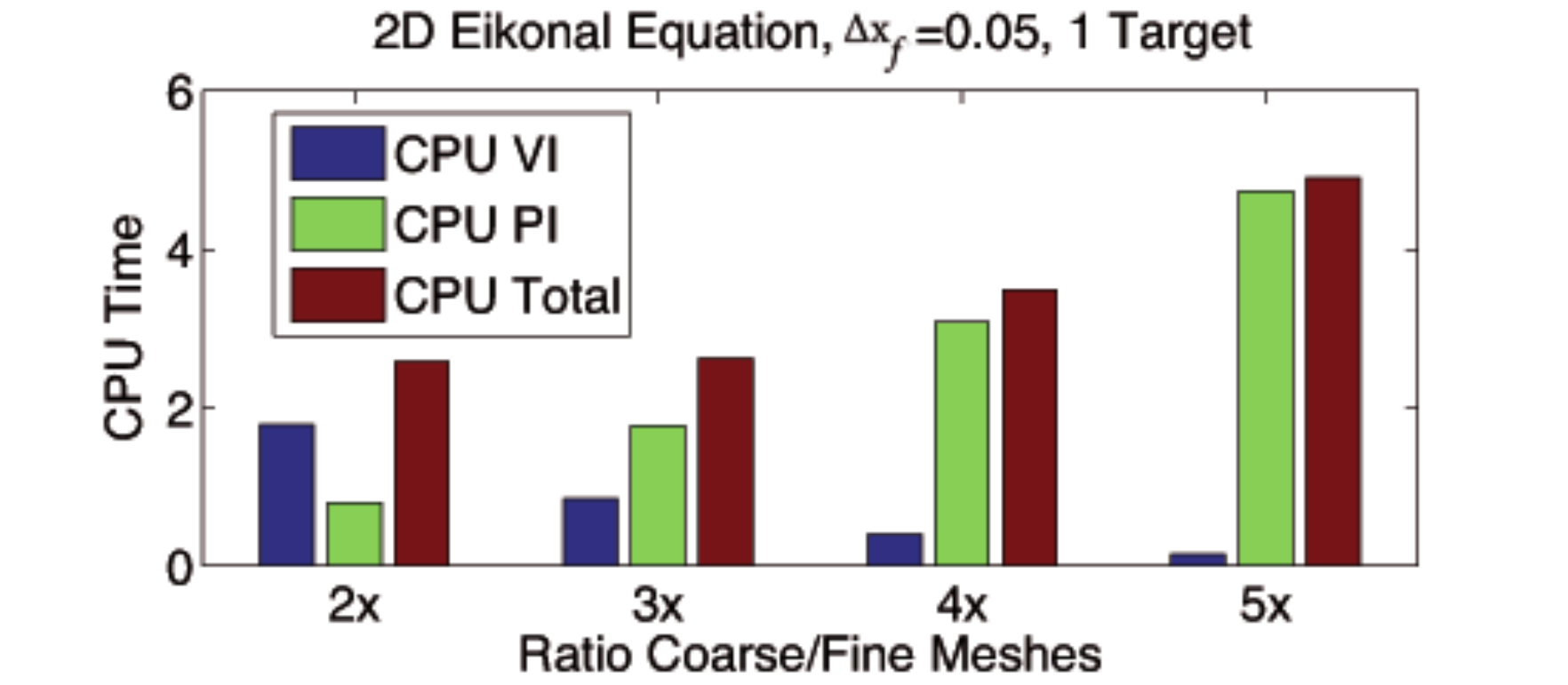,width=0.5\linewidth,clip=} &
\epsfig{file=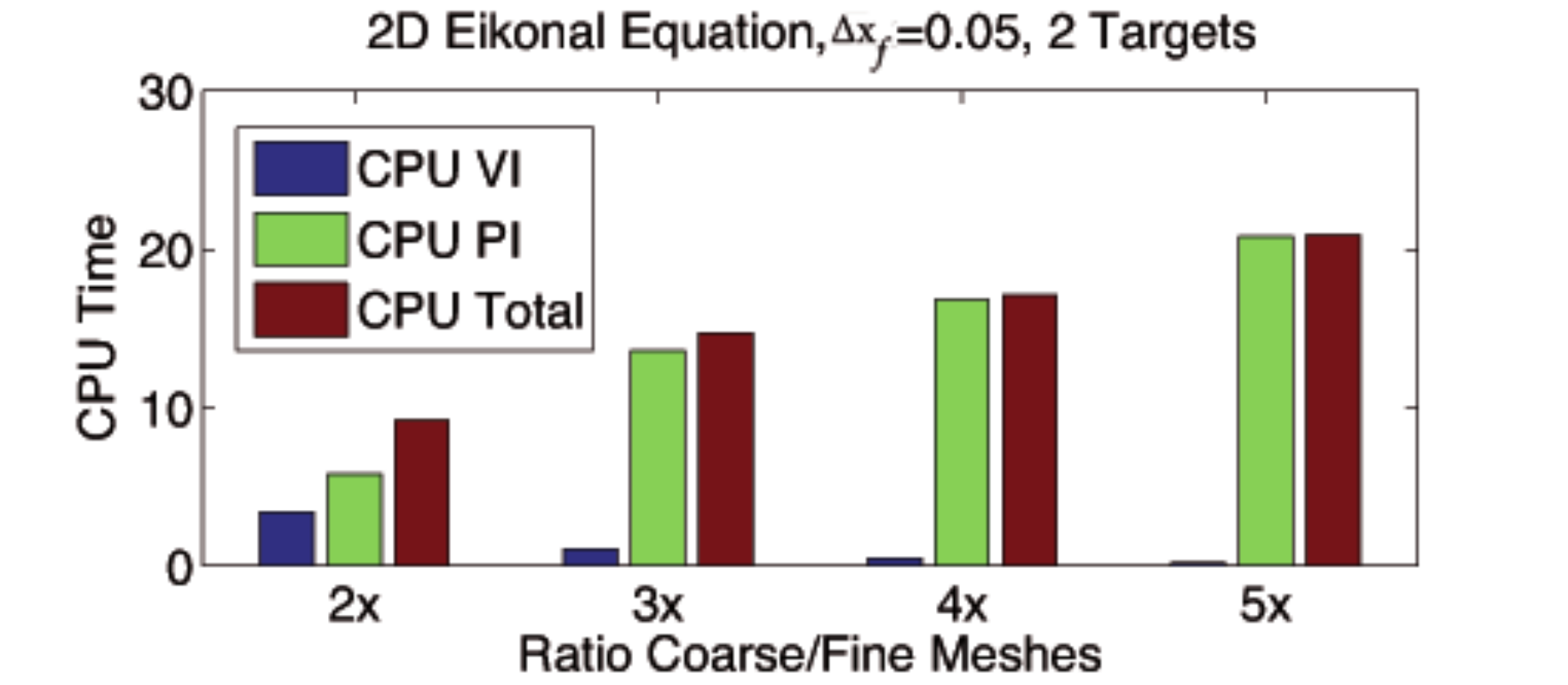,width=0.5\linewidth,clip=} \\
\epsfig{file=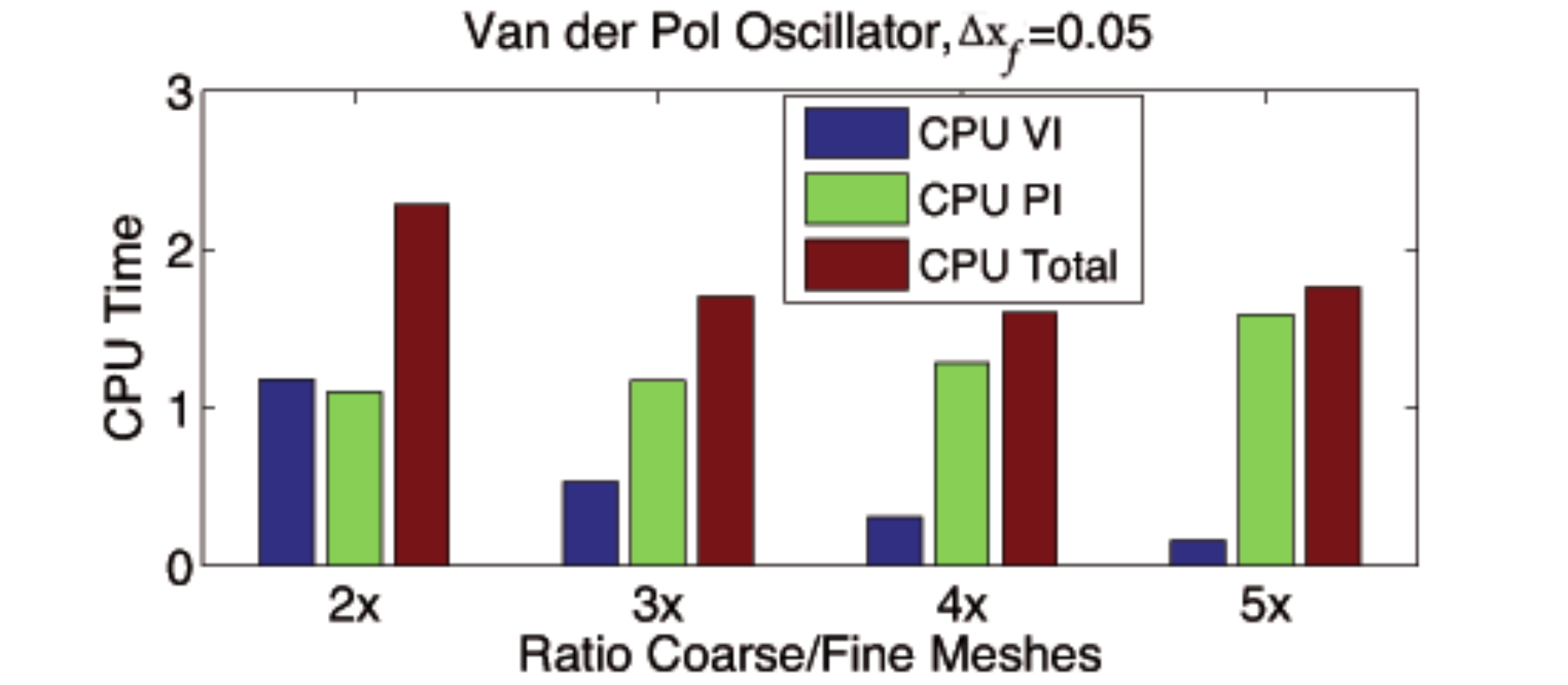,width=0.5\linewidth,clip=} &
\epsfig{file=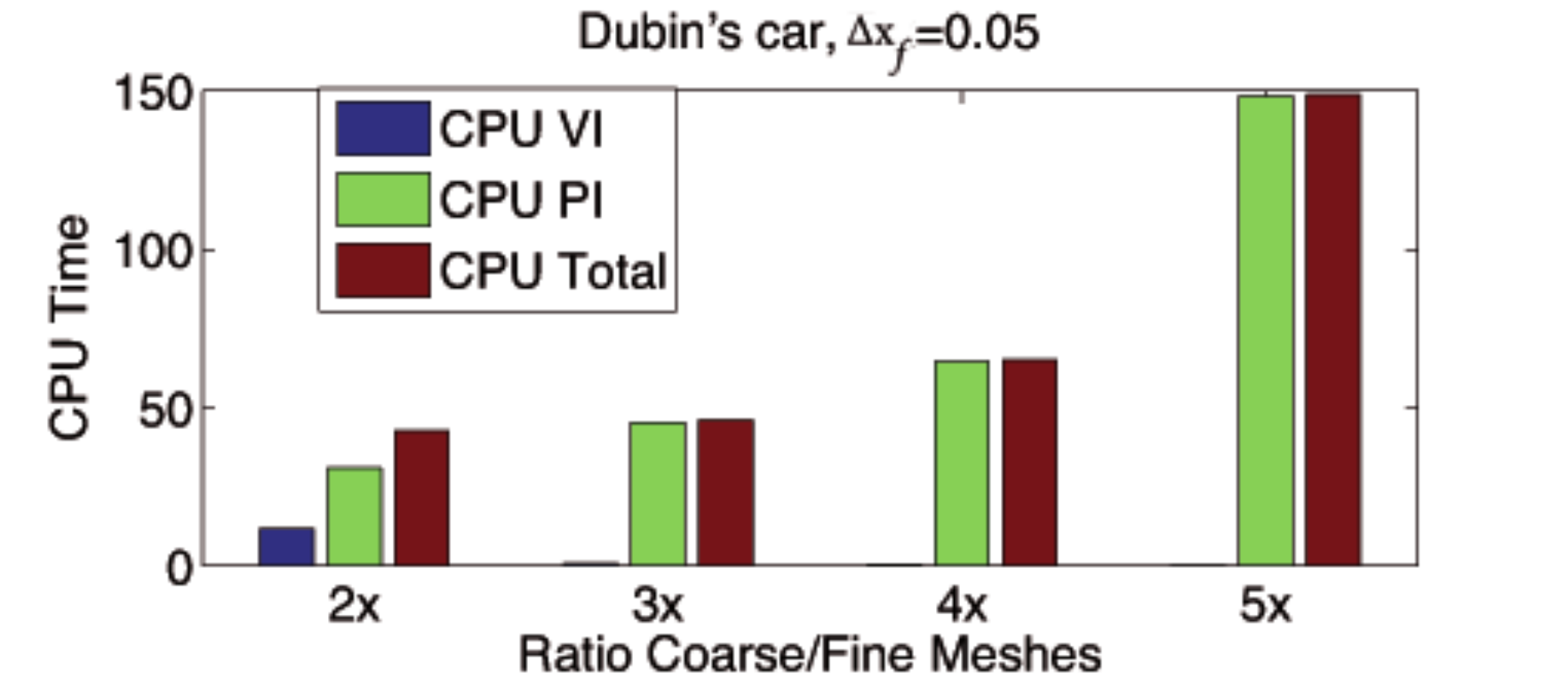,width=0.5\linewidth,clip=}
\end{tabular}}
\caption{Ratios $\Delta x_{c}/\Delta x_f$ and CPU time for different control problems. A good overall balance can be observed in most cases by considering $\Delta x_c=2\Delta x_f$.}
\label{fig:ratios}
\end{figure}

\noindent\textbf{Coarse and fine meshes.} The main trade-off of the accelerated algorithm is related to this point. For a good behavior of the PI part of the algorithm, a good initialization is required, but this should be obtained without deteriorating the overall performance. Too coarse VI will lead to poor initialization, while fine VI will increase the CPU time. We recall that for this paper we assume regular equidistant meshes with mesh parameter $\Delta x$. If we denote by $\Delta x_c$ and by $\Delta x_f$ the mesh parameters associated to the coarse and fine grids respectively,  numerical findings illustrated in Figure \ref{fig:ratios} suggest that for minimum time problems and infinite horizon optimal control, a good balance is achieved with $\Delta x_c=2\Delta x_f$. In the case of minimum time problem, additionally, it is important that the coarse mesh is able to accurately represent the target.

\noindent\textbf{Accuracy.} Both VI and PI algorithms require a stopping criterion for convergence. Following \cite{SR04}, the stopping criteria is given by
\begin{equation*}
||V^{k+1}-V^k||\leq C \Delta x^2\,,
\end{equation*}
which relates the error to the resolution of the state-space mesh. The constant $C$ is set to $C=\frac15$ for the fine mesh, and for values ranging from 1 to 10 in the coarse mesh, as we do not strive for additional accuracy that usually will not improve the initial guess of the control field. However, different options have been extensively discussed in the literature, as in \cite{S97} for instance, where the stopping criteria is related to a variability threshold on the control space.\\

\noindent\textbf{Policy evaluation.} In  every main cycle of the policy iteration algorithm, provided the interpolation operator is linear, as it is in our case, a solution of the linear system \eqref{policy_ev} is required. This can be performed in several ways, specially given the sparsity of the system. For sake of simplicity and in order to make numerical comparisons with the VI scheme, we use a fixed point iteration, i.e., the policy evaluation is implemented as
\begin{equation}
V_i^{k,j+1} = \Delta t g\left(x_i,a_i^k\right)+e^{-\lambda \Delta t} I\left[V^{k,j}\right]\left(x_i+\Delta t f\left(x_i,a_i^k\right)\right)
\end{equation}
with initial guess $V^{k,0}=V^{k-1,\infty}$. We use the same stopping criteria as for the global iteration.

\noindent\textbf{Minimization.} Although counterexamples can be constructed in order to show that it is not possible to establish error bounds of the PI algorithm independently of the (finite) number of controls \cite{SR04}, the algorithm does not change its  performance when the control set is increased, and therefore the \textsl{argmin} computation required for the policy update can be performed by discretizing the set of controls and evaluating all the possible arrival points. Note that, in order to avoid the discretization of the control set, minimizers can be computed using  Brent's algorithm, as in \cite{CFF04}.

\noindent\textbf{A remark on parallelism.} Although the numerical tests that we present were performed in a serial code, we note that the accelerated algorithm allows an easy parallel implementation. Whenever an iterative procedure is performed over the value function, parallelism can be implemented via a domain decomposition of the state space as in \cite{FLS94,CCFP12}. If a control space discretization is also performed, the policy update \eqref{eq:polup} can also be parallelized with respect of the set of controls.

\section{An Error Estimate for API}
Let us present a theoretical estimate for the error related to the API algorithm. Since the method is based on the coupling of a VI approximation with a PI approximation, we need to couple the general error estimate for VI  \cite[Appendix A]{BCD97} with the a local error estimate for the PI method. We will distinguish two cases: the smooth case and nonsmooth case.
In the first case we will use the classical local error estimate for Newton method whereas in the nonsmooth case we will use the error estimate proved by Santos and Rust in \cite{SR04}. Note that the main contribution of Santos and Rust has been to extend the results by Puterman and Brumelle \cite{PB79} to the fully discrete problem where the algorithms are applied on a grid at discrete time values $t_n$. \\
Here  we work on two fully discrete methods and the crucial point is always  to decide when to switch from the first method to the second method guaranteeing convergence.

\noindent Let us start recalling a classical result for Newton's method. More details can be found in \cite{D04}. Assume we have to solve a nonlinear operator equation
$$F(x)=0,$$ where $F:D\rightarrow Y$, $D\subset X$ and $X,Y$ are Banach spaces. Let $F$ be at least once continuously differentiable and let us denote by $F'$ its Fr\'echet derivative. Suppose we have a starting guess $x_0$ of the unknown solutions $x^*$ at hand. Then the successive linearization leads to the general Newton method
\begin{equation}\label{eq:New}
F'(x_k) d_k=-F(x_k), \quad \hbox{where }d_k=x_{k+1}-x_k,\quad k=0,1,\ldots
\end{equation}
so we have to solve a linear system at each iteration to compute the direction of displacement $d_k$. The following theorem provides sufficient conditions to ensure the convergence of the method for a given initial guess $x_0$.
\begin{theorem}\label{th:New}
Let $F:D\rightarrow Y$ be a continuously Fr\'echet differentiable mapping with $D\subset X$ open and convex. For a starting point $x_0\in D$ let $F'(x_0)$ be invertible. Assume that the following conditions hold true,
\begin{equation}\label{new1}
\|F'(x_0)^{-1}F(x_0)\|\le\alpha
\end{equation}
\begin{equation}\label{new2}
\|F'(x_0)^{-1}(F'(y)-F'(x))\|\le \bar{\omega}\|y-x\|\quad x,y\in D
\end{equation}
\begin{equation}\label{new3}
h_0:=\alpha\bar{\omega}\le \frac{1}{2}
\end{equation}
\begin{equation}\label{new4}
B(x_0,\rho)\subset D\qquad \rho:=\dfrac{1-\sqrt{1-2h_0}}{\bar{\omega}}
\end{equation}
Then, the sequence $\{x_k\}$ obtained by \eqref{eq:New} is well-defined, remains in $B(x_0,\rho),$ and converges to some $x^*$ with $F(x^*)=0.$ For $h_0<\frac{1}{2},$ the convergence is quadratic.
\end{theorem}
In the sequel of this section we will use particular notations to avoid confusion. In particular, we will denote by lower case letters the  functions (as the value function which will be indicated by $v$) defined everywhere in the domain and with capital letters vectors which will typically correspond to result of the fully discrete methods. In particular, $V$ will be the fully discrete solution generated by the VI-method and $W$ will be the fully discrete solution generated by the PI-method (note that on a fixed grid $V$ and $W$ coincide). Clearly the cardinality of the vectors will be the cardinality of the grid nodes, i.e. is equal to $card (I_G)$ where $I_G$ denoted the set of indices belonging to the grid.
We will also denote by $\Delta=(\Delta t, \Delta x)$ the vector containing the parameters of the discretization, clearly the vectors $V$ and $W$ will depend on $\Delta$ as well as the the approximating sequences $V^n$, $W^n$  generated by the two methods. To avoid cumbersome notations, we will explicitly stress this dependence only when it is necessary.
\noindent Note that the Policy Iteration scheme may be considered as a Newton-like algorithm (as explained in \cite{PB79}) and a more precise result  for the fully discrete scheme has been proved by Santos and Rust in \cite{SR04} (at  p. 2107) establishing quadratic convergence.

\begin{theorem}\label{PI:quad}
Assume $V$ is the fixed point of the discretized Bellman equation \eqref{hjbhk} (which can be numerically approximated either by VI and PI methods). Assume that the continuous value function $v:\Omega\rightarrow \R$ is concave and let $\{W^n\}_{n\ge1}$ be a sequence of functions generated by the PI-method  \eqref{policy_ev} and \eqref{eq:polup} with local reconstruction based on linear interpolation on the grid.
Moreover, let  the initial guess $W^0$ be close enough to the solution $V$, we have:
\begin{equation}\label{PI:est}
\|V - W^{n+1}\|_\infty \le C\dfrac{\gamma}{\Delta x ^2 (1-\gamma)}\|V- W^n\|_\infty^2.
\end{equation}
where $\gamma$ is a positive parameter in $(0,1)$.
\end{theorem}

\noindent In our case, $\gamma=e^{-\lambda \Delta t}$ is the contraction coefficient of the nonlinear operator $T$.

\noindent The above result guarantees convergence provided the approximate value function is close enough to the solution $V$,  so we to enter into a small neighborhood of $V$ via the VI-method in order to have a good initialization. Let us recall  the error estimate for the value iteration scheme \eqref{hjbhk}.

\noindent We can rewrite relation \eqref{hjbh} with compact notation in a fixed point form $V_i=T(V_i)$, for $i=1,\dots, N_G$
where the map $T:\R^{N_G}\rightarrow\R^{N_G}$ is defined componentwise as
$$\left(T(V)\right)_i\equiv\min_{a\in A} [e^{-\lambda \Delta t} \Lambda(a)V+\Delta t g(a)]_i$$
where $g_i(a)=g(x_i,a),$ and $\Lambda(a)\in\R^{{N_G}\times {N_G}}$ matrix and we have:
$$v(x_i+\Delta t f(x_i, a))=\sum_{j=1}^{N_G} \lambda_{ij}(a)v(x_j),\quad i=1,\ldots,N_G.$$
The following theorem gives an estimate of the error of the value function in the Value Iteration algorithm provided the invariant condition
$x_i+\Delta t f(x_i,a)\in \Omega$, for any $i$,  is satisfied (proof and other details can be found in \cite[Appendix A]{BCD97}):
\begin{theorem}\label{VI:conv}
Let $v$ be the continuous value function and $V$ be the solution of (\ref{hjbhk}). Assume that:
$$
f:\R^n\times A\rightarrow\R^n\hbox{ and } g:\R^n\times A\rightarrow\R\hbox{ are continuous,}
$$
$$
\|f(x,a)-f(y,a)\|\leq L_f \|x-y\| \hbox{ for any }a\in A\hbox{ and } \|f\|_\infty\le M_f,
$$
$$
\|L(x,a)-L(y,a)\|\leq L_g\|x-y\| \hbox{ for any }a\in A\hbox{ and } \|g\|_\infty\le M_L,
$$
where $M_L, M_f, K_L>0.$  Let us assume that $ \lambda>L_f$, the following inequality holds:
\begin{equation}\label{hjb:est}
\max_{i\in {N_G}} \|v(x_i)-V_i\|\le C (\Delta t)^{1/2}+\dfrac{L_f}{\lambda(\lambda-L_f)}\dfrac{\Delta x}{\Delta t}.
\end{equation}
\end{theorem}
\medskip
Since $T$ is a contraction mapping in $\R^{N_G}$ given $V^0$, the sequence
$$V^n=T(V^{n-1}),\quad n=1,2,\ldots$$
will converge to $V$ {\em for every initial condition} $V^0\in\R^d.$ Moreover, under the assumptions of
Theorem \ref{VI:conv} we have the following relation between the current iteration $V^n$ and the initial condition $V^0:$
\begin{equation*}
\|V-V^n\|_\infty \le e^{-n \lambda \Delta t} \|V-V^0\|_\infty = \gamma ^n \|V-V^0\|_\infty
\end{equation*}
\par
\noindent Let us make a rather important remark. Let us suppose we have found on a grid with uniform space step $\Delta x$  an approximate value function such that
\begin{equation}\label{VI_end}
\|V^{n+1,\Delta x}-V^{n,\Delta x}\|_\infty\le\varepsilon\,.
\end{equation}

\noindent We need a result regarding the local reconstructions based on linear interpolation since we typically use two grid, one coarse and one fine grid and we project the $V$ computed on the coarse grid onto the fine grid. Now let us suppose that we have for the fine grid the nodes are $G=(x_1\ldots,x_{i-1},x_{i-1/2},x_i,x_{i+1/2},x_{i+1},\ldots x_{N_G} )$  and that $\Delta x=x_i-x_{i-1}$ whereas $x_{i}-x_{i-1/2}=\frac{\Delta x}{2}$ (i.e. the fine grid is obtained just adding the middle points of every cell).
First of all we want to understand what happens if we initialize the PI method with
$W^0=V^{n+1,\Delta x/2}$ which is obtained by linear interpolation from $V^{n+1,\Delta x}.$ In the node of the grid with step size $\Delta x$ we will find the same value function, e.g $V^{n+1,\Delta x/2}_i=V^{n+1,\Delta x}_i$
 whereas in the new points of the grid we get $V^{n,\Delta x/2}_{i\pm1/2}=({V^{n,\Delta x}_i+V^{n,\Delta x}_{i\pm 1}})/2$.
Then we can obtain the same accuracy of \eqref{VI_end}, for example for $i+1/2$ we get:
\begin{eqnarray*}
\|V^{n+1,\Delta x/2}_{i+1/2}-V^{n,\Delta x/2}_{i+1/2}\|_\infty&=&\left\| \dfrac{V^{n+1,\Delta x}_i+V^{n+1,\Delta x}_{i+1}}{2}-\left(\dfrac{V^{n,\Delta x}_i+V^{n,\Delta x}_{i+1}}{2}\right)\right\|_\infty\\
&=&\dfrac{1}{2}\left\|\left( V^{n+1,\Delta x}_{i+1}-V^{n,\Delta x}_{i+1}\right)+\left(V^{n+1,\Delta x}_i-V^{n,\Delta x}_i\right)\right\|_\infty\\
&\le&\dfrac{1}{2}\left(\left\|\left( V^{n+1,\Delta x}_{i+1}-V^{n,\Delta x}_{i+1}\right)\|_\infty+\|\left(V^{n+1,\Delta x}_i-V^{n,\Delta x}_i\right)\right\|_\infty\right)\\
&\le&\varepsilon
\end{eqnarray*}
In this way we have shown that the Value function obtained by linear interpolation from the coarse grid has the same error and we do not loose accuracy restarting the procedure for the PI-method on the fine grid. The controls are simply obtained computing the \textsl{argmin} of the new value function. Of course, we will have the same controls in $x_i$ then we need to compute half of them so this procedure  is not expensive.
$$a_{i+1/2}=\arg\min_{a\in A}\left\{e^{-\lambda\Delta t}V^{n+1,\Delta x/2}(x_{i+1/2}+\Delta  t\,f(x_{i+1/2},a)+\Delta  t L(x_{i+1/2},a)\right\}$$

\noindent Finally, let us combine the above result for the value (Theorem \ref{VI:conv})  and the policy iteration method (Theorem \ref{PI:quad})  in order to get a global convergence result for the API algorithm. For simplicity we drop the subscript $\Delta x$ over the approximate values $V$ and $W$.\\
We call $V^k_{API}$ the current iteration of the accelerated scheme, which is obtained combining $N$ iteration of the the VI-method with $V^0$ as  initial
 condition and  $M$ iterations for the PI-method restarted from $W^0=V^{N, \Delta x/2}$. Clearly $V^k_{API} =W^M$.

\begin{theorem}
Let all the assumptions for Theorem \ref{th:New}, Theorem \ref{PI:quad} and Theorem \ref{VI:conv} hold true. We have the following result:\\
1) For any given positive $\rho$, there exists an $N$ such that $\max_{i\in I_G} |v(x_i)-V^N|< \rho$;\\
2) The sequence generated by the PI-method, with initial condition $W^0\equiv V^N$,  will converge quadratically to $V$ with the error estimate
\begin{equation}\label{API-estimate}
\|V^k_{API}-V\|\le  \mu \left(\prod _{k=1}^M \mu^{2(k-1)}\rho^2\right ).
\end{equation}
where
\begin{equation}
\mu =C \dfrac{ \gamma}{\Delta x^2 (1-\gamma)}.
\end{equation}
\end{theorem}
{\em Proof}
1) Let  $V^0$ be the initial condition for the VI-procedure, then by a contraction argument we have
\[\|V^N -V\|\le e^{-N \lambda \Delta t}\|V^0-V\|=\gamma ^N \|V^0-V\|\]
Moreover, by applying \eqref{VI:conv}, we can get the explicit estimate  \eqref{hjb:est}  with respect to the discretization parameters $\Delta$ so we can guarantee that we enter any small ball of radius $\rho$ centered at the exact value function.
This allows to start up the PI-method in the convergence region.\\
2)  To get our estimate let us start with the \eqref{PI:est}. This implies, in the quadratic convergence region,
\begin{equation}
\Vert V-W^M \Vert _\infty \le\mu \|V-W^{M-1}\|_\infty ^2\le \dots \le  \mu \left(\prod _{k=1}^M \mu^{2(k-1)}\right ) \Vert V-W^0\Vert_\infty ^ {2M}
\end {equation}
Clearly, if $ \Vert V-W^0\Vert_\infty \le \rho< 1$ we get
\begin{equation}
\Vert V-W^M \Vert _\infty \le\mu \left(\prod _{k=1}^M \mu^{2(k-1)}\right) \rho^{2M}=\mu  \prod _{k=1}^M( \mu^{k-1}\rho)^2
\end {equation}
Note that the upper bound on the estimate increases as far as $\mu$ increases but the condition $\mu^{M-1}\rho<1$ is sufficient to guarantee convergence.

\vspace{0.2cm}\noindent
{\em A remark on the (ideal) smooth case.}\\
Note that in general we are computing non smooth solution and the operator  $F(U):=U-T(U)$ is not regular because we are using local linear interpolation on the grid. However,  applying other local interpolation operators and assuming that the feedback control is smooth  in a neighborhood of the solution we could directly apply  the Newton method and again the initial condition has a crucial role. We can now couple the general result for Newton with the error estimate for the VI-method.

\noindent Therefore, we need to check the hypotheses of Theorem \eqref{th:New} are satisfied for $F(V)=0$. Newton method writes as
$$J_F(W^n)(W^{n+1}-W^n)=-F(W^n)\,,$$
where the Jacobian of F is $J_F(U):=I-J_T(U),$ and $I$ is the identity matrix.

\noindent The residual $r^n$ and the error $e^n$ at the $n-$th iteration are given by:
$$r^n\equiv V^n-T(V^n)=V^n-V^{n-1}\qquad e^n\equiv \|V^n-V\|.$$
Then, we fix $\eta>0$ such that:
$$\|r^N\|\le\eta \hbox{  if and only if } \|V^N-V\|\le \rho.$$
An explicit expression for $\rho$ is given in in Theorem \eqref{th:New}. Then, recalling that we initialize $W^0=V^N$ it is required that
$$\|J_F(V^N)^{-1} F(V^N)\|\le \alpha\,,$$
which is guaranteed since $V$ is the solution coming from the value iteration scheme which is known to be convergent for any initial condition.
To guarantee convergence, we also have to assume the following inequality:
$$\|J_F(V^N)^{-1}(J_F(W)-J_F(U))\|\le \bar{\omega}\|W-U\|\,,\quad\forall \;W,U\in B(V^N,\rho).$$


\section{Numerical tests}
This section presents a comprehensive set of tests assessing the performance of the proposed accelerated algorithm. We compare the results with solutions given by the classical value iteration algorithm, policy iteration, and the accelerated monotone value iteration method. In some examples we also include an accelerated algorithm based on a monotone value iteration in the set of subsolutions (AMVI), as presented in \cite[Appendix A]{BCD97}, and a Gauss-Seidel variation of this method (GSVI) as in \cite{G96}. In a first part we develop tests related to infinite horizon optimal control, to then switch to the study of minimum time problems. We conclude with an extension to applications related to optimal control of partial differential equations. We focus on grid resolution, size of the discretized control space, performance in presence of linear/nonlinear dynamics, targets, and state space dimension. All the numerical simulations reported in this paper have been made on a MacBook Pro with 1 CPU Intel Core i5 2.3 Ghz and 8GB RAM.\\

\noindent\textbf{Infinite horizon optimal control problems}
\subsection{Test 1: A non-smooth 1D value function}
\begin{figure}[ht]
\begin{center}
\includegraphics[scale=0.4]{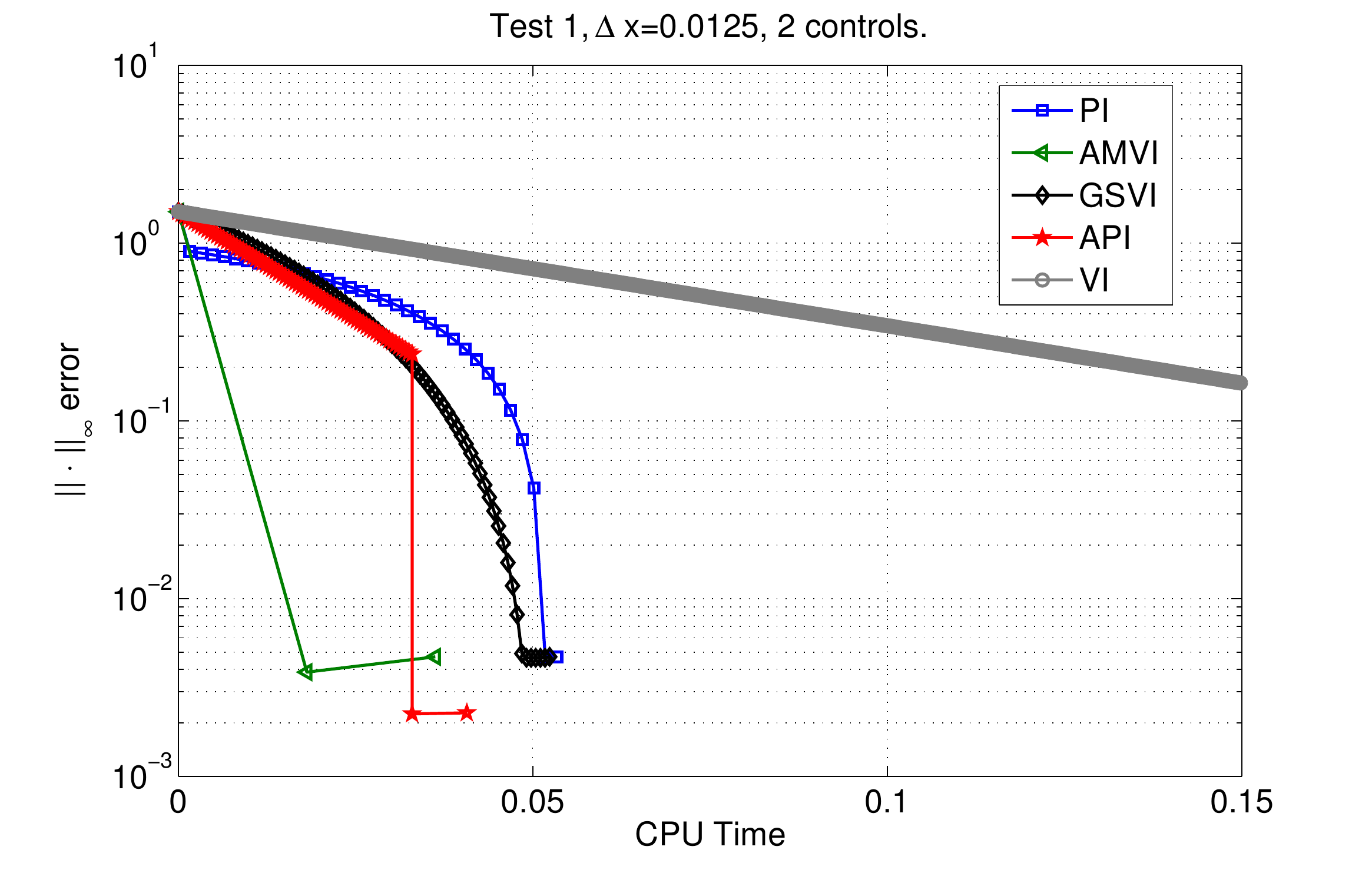}
\end{center}
\caption{Test 1 (non-smooth value function): error evolution for different algorithms.}
\label{fig:Test1}
\end{figure}

We first consider a one-dimensional optimal control problem appearing in \cite[Appendix A]{BCD97}. Using a similar notation as in Section \ref{sec2}, we set the computational domain $\Omega=]-1\,,1[$, the control space $A=[-1\,,1]$, the discount factor $\lambda=1$, the system dynamics $f(x,a)=a(1-|x|)$, and the cost function $g(x,a)=3(1-|x|)$. The exact optimal solution for this problem is
\begin{align*}
v(x)=\left\{\begin{array}{cc}
\frac32(x+1)&\text{for }\, x<0\,,\\[1.5ex]
\frac32(1-x)&\text{elsewhere}\,,
\end{array}\right.
\end{align*}
which has a kink at $x=0$. We implement every proposed algorithm, and results concerning CPU time and number of iterations are shown in Table \ref{tab:Test1}; for different mesh configurations, we set $\Delta t=.5\Delta x$ and we discretize the control space into a set of 20 equidistant points. The notation VI$(2 \Delta x)$ in Table \ref{tab:Test1} stands for the computation of the solution with a VI method considering a coarse grid of $2\Delta x.$ Then it is applied the PI method with a stepsize $\Delta x.$ (PI($\Delta x$) in the table). This notation, in the table, is kept in all the tests. In this test case, as expected, we observe that the VI algorithm is always the slowest option, with iteration count depending on the number of mesh nodes; this feature is also observed for the PI algorithm, although the number of iterations and CPU time are considerably smaller. On the other hand, the AMVI scheme has an iteration count independent of the degrees of freedom of the system, with an almost fixed CPU time, as the time spent on fixed point iterations is negligible compared to the search of the optimal update direction. In this particular example, the exact boundary conditions of the problem are known ($v(x)=0$ at $\partial\Omega$) and it is possible to construct monotone iterations by starting from the initial guess $v(x)=0$.
The GSVI method exhibits a similar performance as the PI algorithm, with a considerably reduced number of iterations when compared to VI. Note however, that this implementation requires the pre-computation and storage of the interpolation coefficients and a sequential running along the mesh, which can be unpractical of high-dimensional problems are considered. Finally, the API algorithm exhibits comparable CPU times as AMVI, performing always better than VI, PI and GSVI. In this particular case, the choice of the mesh ratio between the coarse and fine meshes can be suboptimal, as the time spent on the VI coarse pre-processing represents an important part of the overall CPU time. More details on the error evolution throughout the iterations can be observed  in Figure \ref{fig:Test1}; note that the error evolution is measured with respect to the exact solution and not with respect to the next iteration. This latter figure illustrates, for both problems, the way in which the API idea acts: pre-processing of the initial guess of PI leads to proximity to a quadratic convergence neighborhood, where this algorithm can converge in a reduced number of iterations; the fast error decay that coarse mesh VI has in comparison with the fine mesh VI is clearly noticeable.
\noindent In Table \ref{tab:Test1lambda}, we show the performance evolution of the different algorithms when the parameter $\lambda$ decreases to zero. It is expected that for methods based on a fixed point iteration of the value function, the number of iterations required to reach a prescribed error level will gradually increase. This is clearly observed for VI, PI and API, whereas AMVI and GSVI are able to circumvent this difficulty, leading to a constant number of iterations independent of the parameter $\lambda$. Nevertheless, in the overall cpu time, GSVI and API exhibit a similar asymptotic performance.

\begin{table}[ht]
\begin{center}
\ra{1.2}
\resizebox{\textwidth}{!}{
\begin{tabular}{ccccccccccccccc}
\toprule
$\#$ nodes & $\Delta x$ & &VI & &PI & &AMVI & &GSVI & & VI($2\Delta x$)& PI($\Delta x$) & API \\
\cmidrule{1-1}\cmidrule{2-2}\cmidrule{4-4}\cmidrule{6-6}\cmidrule {8-8} \cmidrule{10-10} \cmidrule{12-14}
81  &2.5E-2      & & 9.88E-2 (228)&   & 2.02E-2 (10)&  & 1.99E-2 (3)& &2.25E-2 (41) & &  5.31E-3(23) & 5.22E-3 (2)& 1.05E-2\\
161 &1.25E-2     & & 0.41 (512)   &   & 5.88E-2 (34)&  & 3.8E-2 (3) & &7.71E-2 (81) & &  3.21E-2(73) & 1.73E-2 (2)& 4.94E-2\\
321 &6.25E-3     & & 1.89 (1134)  &   & 0.21 (65)   &  & 7.48E-2 (3)& &0.29 (161)   & &  0.16(200)   & 2.62E-2 (2)& 0.19\\
\bottomrule
\end{tabular}}
\vskip 2mm
\caption{Test 1 (1D non-smooth value function): CPU time (iterations) for different algorithms.}\label{tab:Test1}
\end{center}
\end{table}

\begin{table}[ht]
\begin{center}
\ra{1.2}
\resizebox{\textwidth}{!}{
\begin{tabular}{ccccccccccccc}
\toprule
$\lambda$ & &VI & &PI & &AMVI & &GSVI & & VI($2\Delta x$)& PI($\Delta x$) & API \\
\cmidrule{1-1}\cmidrule{3-3}\cmidrule{5-5}\cmidrule {7-7} \cmidrule{9-9} \cmidrule{11-13}
1    & & 1.31 (1134)&   & 0.16 (65) &  & 5.73E-2 (3)& &0.19 (161) & &  7.41E-2 (112) & 3.20E-2  (2)& 0.11\\
0.1  & & 2.45 (2061)&   & 0.46 (138)&  & 5.82E-2 (3)& &0.19 (161) & &  0.12 (203)    & 5.41E-2  (2)& 0.18\\
1E-2 & & 2.63 (2244)&   & 0.67 (159)&  & 6.18E-2 (3)& &0.19 (161) & &  0.12 (220)    & 6.483E-2 (2)& 0.19\\
1E-3 & & 2.65 (2265)&   & 0.74 (161)&  & 7.75E-2 (4)& &0.19 (161) & &  0.13 (222)    & 6.41E-2  (2)& 0.19\\
\bottomrule
\end{tabular}}
\vskip 2mm
\caption{Test 1 (1D non-smooth value function): CPU time (iterations) for different algorithms and different values of $\lambda$, in a fixed mesh with $321^2$ nodes and 2 control values.}\label{tab:Test1lambda}
\end{center}
\end{table}

\subsection{Test 2: Van Der Pol oscillator} In a next step we consider two-dimensional, nonlinear system dynamics given by the Van der Pol oscillator:
\begin{equation*}
f(x,y,a)=\left(\begin{array}{c}
y\\
(1-x^2)y-x+a\end{array}\right)\,.
\end{equation*}
Remaining system parameters are set:
\begin{equation*}
\Omega=]-2\,,2[^2\,,\quad A=[-1\,,1]\,,\quad\lambda=1\,,\quad \Delta t=0.3\Delta x\,,\quad g(x,y,a)=x^2+y^2\,,
\end{equation*}
and the control space is discretized into 32 equidistant points. We perform a similar numerical study as in the previous example, and results are shown in Table \ref{tab:Test2}. For computations requiring an exact solution, we consider as a reference a fine grid simulation with $\Delta x= 0.00625$.

\noindent We set a constant boundary value $v(x)=3.5$ at $\partial \Omega$, which can be interpreted as a penalization on the state. If accurate solutions near the boundary are required, a natural choice in our setting would be to perform simulations over an enlarged domain and then restrict the numerical results to a subset of interest.
From this test we observe a serious limitation on the AMVI algorithm. The number of iterations now depends on the number of nodes, and even though the number of iterations is still lower than in the VI algorithm, the CPU time increases as for every iteration a search procedure is required. As it is not possible to find monotone update directions, the AMVI algorithm becomes a VI method plus an expensive search procedure. This lack of possible monotone update can be due to several factors: the nonlinear dynamics, the existence of trajectories exiting the computational domain, and a sensitivity to the artificial boundary condition. We report having performed similar tests for the linear double integrator problem ($\ddot x=a$) with similar results, therefore we conjecture that in this case, the underperformance of the AMVI scheme is due to poor boundary resolution and its use by optimal trajectories. Unfortunately, this is a recurrent problem in the context of optimal control. This situation does not constitute a problem for the API algorithm, where a substantial speedup is seen in both coarse and fine meshes. Note that compared to PI, the accelerated scheme has a number of iterations on its second part which is independent of the mesh parameters as we are in a close neighborhood of the optimal solution.

\begin{table}[ht]
\begin{center}
\ra{1.2}
\resizebox{\textwidth}{!}{
\begin{tabular}{cccccccccccc}
\toprule
$\#$ nodes & $\Delta x$& &VI & &PI & &AMVI & & VI($2\Delta x$)& PI($\Delta x$) & API \\
\cmidrule{1-1}\cmidrule{2-2}\cmidrule{4-4}\cmidrule {6-6} \cmidrule{8-8}\cmidrule{10-12}
$81^2$ & 5E-2        & & 39.6 (529)&   & 5.35 (8)&  & 1.42E2 (3)& &  1.86 (207)   & 1.47 (4)& 3.33 (211)\\
$161^2$ & 2.5E-2       & & 3.22E2 (1267)&   & 34.5 (11)&  & 1.01E3 (563)& &  10.7(165)  & 6.87 (4)& 17.5 (169)\\
$321^2$  & 1.25E-2      & & 3.36E4 (2892)& & 3.36E2 (14)&  & 1.55E4 (2247)& &  88.9 (451) & 47.7 (4)& 1.36E2 (455)\\
\bottomrule
\end{tabular}}
\vskip 2mm
\caption{Test 2 (Van der Pol oscillator): CPU time (iterations) for different algorithms.}\label{tab:Test2}
\end{center}
\end{table}

\subsection{Test 3: Dubin's Car}
\begin{figure}[ht]
\centering
\resizebox{\textwidth}{!}{
\begin{tabular}{cc}
\epsfig{file=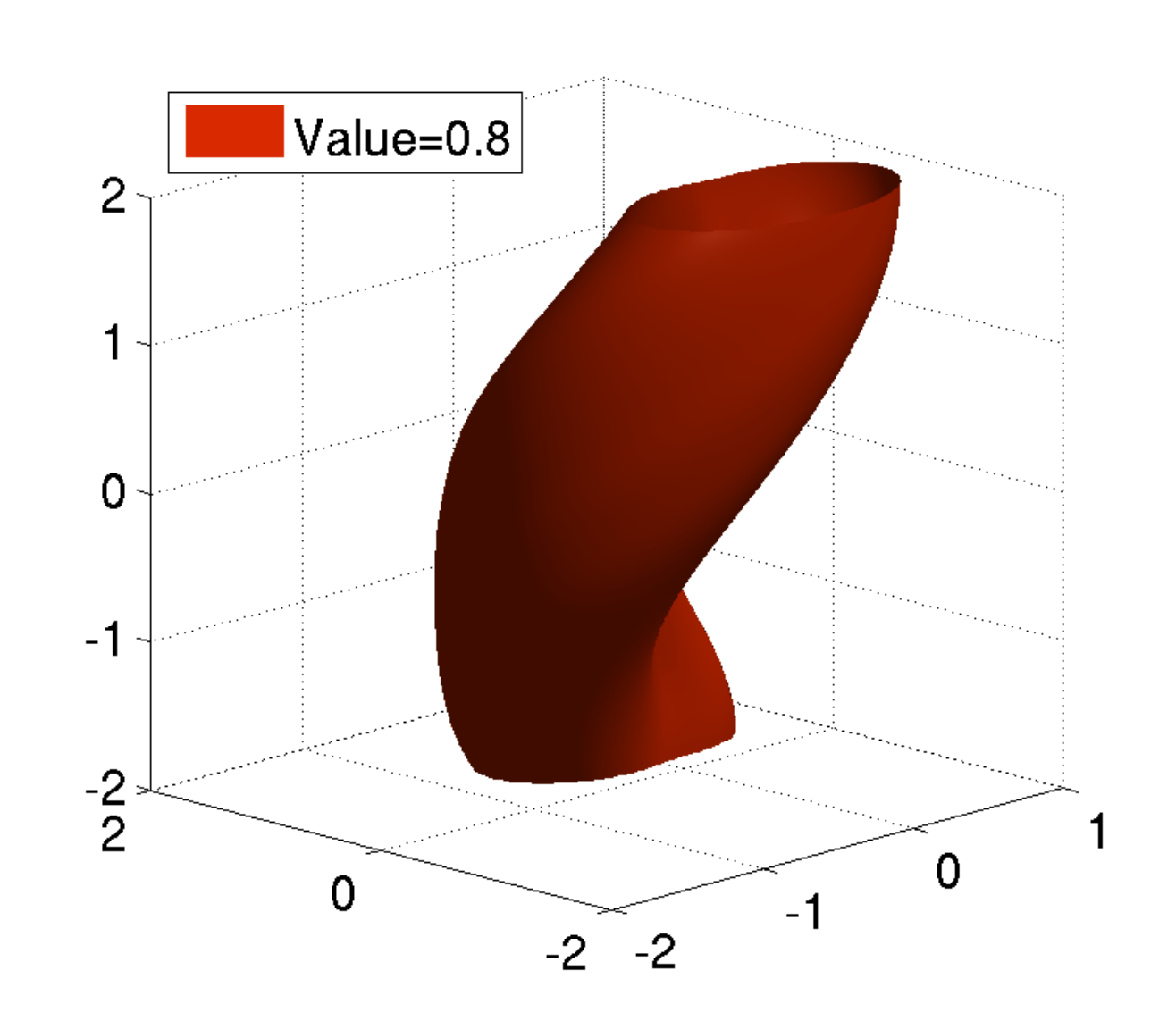,width=0.5\linewidth,clip=} &
\epsfig{file=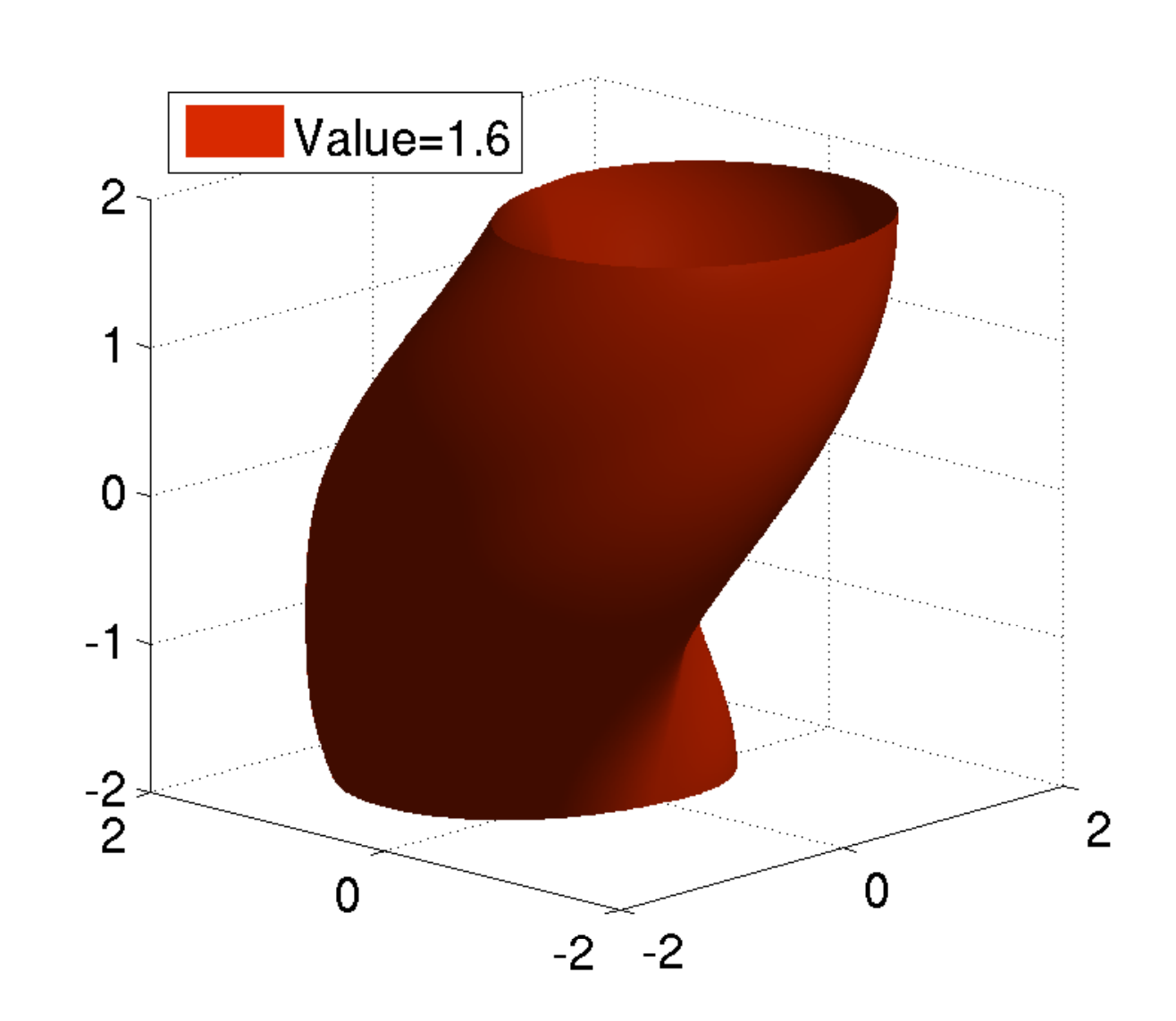,width=0.5\linewidth,clip=}
\end{tabular}}
\caption{Test 3: Dubin's car  value function isosurfaces.}
\label{fig:isodubin}
\end{figure}

Having tested some basic features of the proposed schemes, we proceed with our numerical study of the API method by considering a three-dimensional nonlinear dynamical system given by
\begin{equation*}
f(x,y,z,a)=\left(\begin{array}{c}
\cos(z)\\
\sin(z)\\
a\end{array}\right)\,,
\end{equation*}
corresponding to a simplified version of the so-called Dubin's car, a test problem extensively used in the context of reachable sets and differential games. System parameters are set:
\begin{equation*}
\Omega=]-2\,,2[^2\,,\quad A=[-1\,,1]\,,\quad\lambda=1\,,\quad \Delta t=0.2\Delta x\,,\quad g(x,y,z,a)=x^2+y^2\,,
\end{equation*}
 and the control space is discretized into 11 equidistant points; the boundary value is set to $v(x)=3$ in $\partial \Omega$ and reference solution is taken with $\Delta x=0.0125$. Different isosurfaces for this optimal control problem can be seen in Figure \ref{fig:isodubin}, and CPU times for different meshes are shown in Table \ref{tab:dubin}. This case is an example in which the mesh ratio between coarse and fine meshes is well-balanced, and the time spent in pre-processing via VI is not relevant in the overall API CPU time, despite leading to a considerable speedup of the order of $8\times$ with an order of $10^6$ grid points. In the last line of Table \ref{tab:dubin}, the VI algorithm was stopped after 4 hours of simulation without achieving convergence, which is illustrative of the fact that acceleration techniques in such problems are not only desirable but necessary in order to obtain results with acceptable levels of accuracy.

\begin{table}[ht]
\begin{center}
\ra{1.2}
\scalebox{0.8}{
\begin{tabular}{cccccccccc}
\toprule
$\#$ nodes & $\Delta x$& &VI & &PI & & VI($2\Delta x$)& PI($\Delta x$) & API \\
\cmidrule{1-2}\cmidrule{4-4}\cmidrule{6-6}\cmidrule{8-10}
$41^3$      & 0.1       & & 50.6 (192)&   & 12.2 (12)  & &  0.84 (8)   & 8.52 (3)& 9.36 (11)\\
$81^3$      & 5E-2      & & 1.19E3 (471)&   & 3.28E2 (18) & & 8.98 (39)  & 1.39E2 (9)&1.48E2 (48)\\
$161^3$     & 2.5E-2    & & $\geq$ 1.44E4&   & 9.93E3 (12)  & &  3.02E2 (30)   & 2.92E3 (10)& 2.62E3 (40)\\
\bottomrule
\end{tabular}}
\vskip 2mm
\caption{Test 3 (Dubin's car): CPU time (iterations) for different algorithms}\label{tab:dubin}
\end{center}
\end{table}

\noindent\textbf{Minimum time problems}

\subsection{Tests 4 and 5: minimum time problems in 2D}
The next two cases are based on a two-dimensional eikonal equation. For both problems, common seetings are given by
\begin{equation*}
f(x,y,a)=\left(\begin{array}{c}
\cos(a)\\
\sin(a)\end{array}\right)\,,\quad
A=[-\pi,\,\pi]\,,\quad \Delta t=0.8\Delta x\,.
\end{equation*}
What differentiates the problems is the domain and target definitions; Test 4 considers a domain $\Omega=]-1,\,1[^2$ and a target $\mathcal{T}=(0,0)$, while for Test 5, $\Omega=]-2,\,2[^2$ and $\mathcal{T}=\{x\in\mathbf{R}^2\,:\,||x||_2\leq1\}$. Reference solutions are considered to be the distance function to the respective targets, which is an accurate approximation provided that the number of possible control directions is large enough. For Test 4, with  a discretization of the control space into set of 64 equidistant points, CPU time results are presented in Table \ref{tab:Test4}; it can be seen that API provides a speedup of $8\times$ with respect to VI over fine meshes despite the large set of discrete control points.
\noindent Table \ref{tab:rates} shows experimental convergence rates achieved by the fully discrete scheme, in both $L^1$ and $L^\infty$ norms, which are in accordance with the theoretically expected rate of 1/2.
Test 5 features an enlarged target, and differences in terms of CPU times are presented in Table \ref{tab:Test5} where, for a discrete set of 72 equidistant controls, the speedup is reduced to $4\times$. In general, from a mesh node, larger or more complicated targets represent a difficulty in terms of the choice of the minimizing control, which translates into a larger number of iterations. In this case, the CPU time spent in the pre-processing is significant to the overall CPU time, but increasing this ratio in order to reduce its share will lead to an underperformant PI part of the algorithm.


\begin{table}[ht]
\begin{center}
\ra{1.2}
\resizebox{.8\textwidth}{!}{
\begin{tabular}{cccccccccc}
\toprule
$\#$ nodes & $\Delta x$& &VI & &PI & & VI($2\Delta x$)& PI($\Delta x$) & API \\
\cmidrule{4-4}\cmidrule{6-6}\cmidrule{1-2}\cmidrule{8-10}
$41^2$       & 5E-2        & & 3.16 (37)&   & 1.89 (12)  & &  0.39 (5)   & 0.38 (2)& 0.77 (7)\\
$81^2$      & 2.5E-2       & & 8.23 (69)&   & 4.43 (19) & &  0.80 (12)  & 0.53 (2)& 1.33 (14)\\
$161^2$      & 1.25E-2      & & 39.2 (133)& & 12.6 (13)&   &2.55 (31)  &2.11 (3)  & 4.66 (34)\\
\bottomrule
\end{tabular}}
\vskip 2mm
\caption{Test 4 (2D eikonal): CPU time (iterations) for different algorithms.}\label{tab:Test4}
\end{center}
\end{table}

\begin{table}[ht]
\begin{center}
\ra{1.2}
\resizebox{.7\textwidth}{!}{
\begin{tabular}{ccccccccccc}
\toprule
$\#$ nodes & $\Delta x$ & & $L^1- error$& rate& & $L^\infty- error$ & rate \\
\cmidrule{1-2}\cmidrule{4-5}\cmidrule{7-8}
$41^2$       & 5E-2        & & 2.1E-2  & 0.60& & 8.9E-3 & 0.61\\
$81^2$       & 2.5E-2      & & 1.4E-2  & 0.64& & 5.8E-3 & 0.64\\
$161^2$      & 1.25E-2     & & 8.5E-3  & 0.68& & 3.7E-3 & 0.75\\
$321^2$      & 6.25E-3     & & 5.3E-3  &     & & 2.2E-3 & \\
\bottomrule
\end{tabular}}
\vskip 2mm
\caption{Test 4 (2D Eikonal): Rate of convergence for the API scheme with 64 controls.}\label{tab:rates}
\end{center}
\end{table}

\begin{table}[ht]
\begin{center}
\ra{1.2}
\resizebox{\textwidth}{!}{
\begin{tabular}{cccccccccc}
\toprule
$\#$ nodes & $\Delta x$& &VI & &PI & & VI($2\Delta x$)& PI($\Delta x$) & API \\
\cmidrule{4-4}\cmidrule{6-6}\cmidrule{1-2}\cmidrule{8-10}
$64^2$       & 6.35E-2 & & 4.02 (36)   & & 1.42 (9) & & 0.84 (10)   & 0.53 (4)& 1.37 (14)\\
$128^2$      & 3.15E-2 & & 16.9 (70)   & & 6.25 (14)& & 2.80 (25)  & 1.66 (2)& 4.46 (27)\\
$256^2$      & 1.57E-2 & & 1.09E2 (135)& & 38.7 (16)& & 15.8 (62)  &11.7 (8)  & 27.5 (70)\\
$512^2$      & 7.8E-3  & & 9.80E2 (262)& & 3.98E2(168)& & 1.07E2 (126) & 1.09E2 (12)& 2.16E2 (138)\\
\bottomrule
\end{tabular}}
\vskip 2mm
\caption{Test 5 (2D Eikonal): CPU time (iterations) for different algorithms with 72 controls.}\label{tab:Test5}
\end{center}
\end{table}



\subsection{Tests 6 and 7: minimum time problems in 3D}
We develop a three-dimensional extension of the previously presented examples. System dynamics and common parameters are given by
\begin{equation*}
f(x,y,z,(a_1,a_2))=\left(\begin{array}{c}
\sin(a_1)\cos(a_2)\\
\sin(a_1)\sin(a_2)\\
\cos(a_1)\end{array}\right)\,,\quad
A=[-\pi,\,\pi]\times[0\,,\pi]\,,\quad \Delta t=0.8\Delta x\,.
\end{equation*}

\begin{figure}[ht]
\centering
\resizebox{\textwidth}{!}{
\begin{tabular}{cc}
\epsfig{file=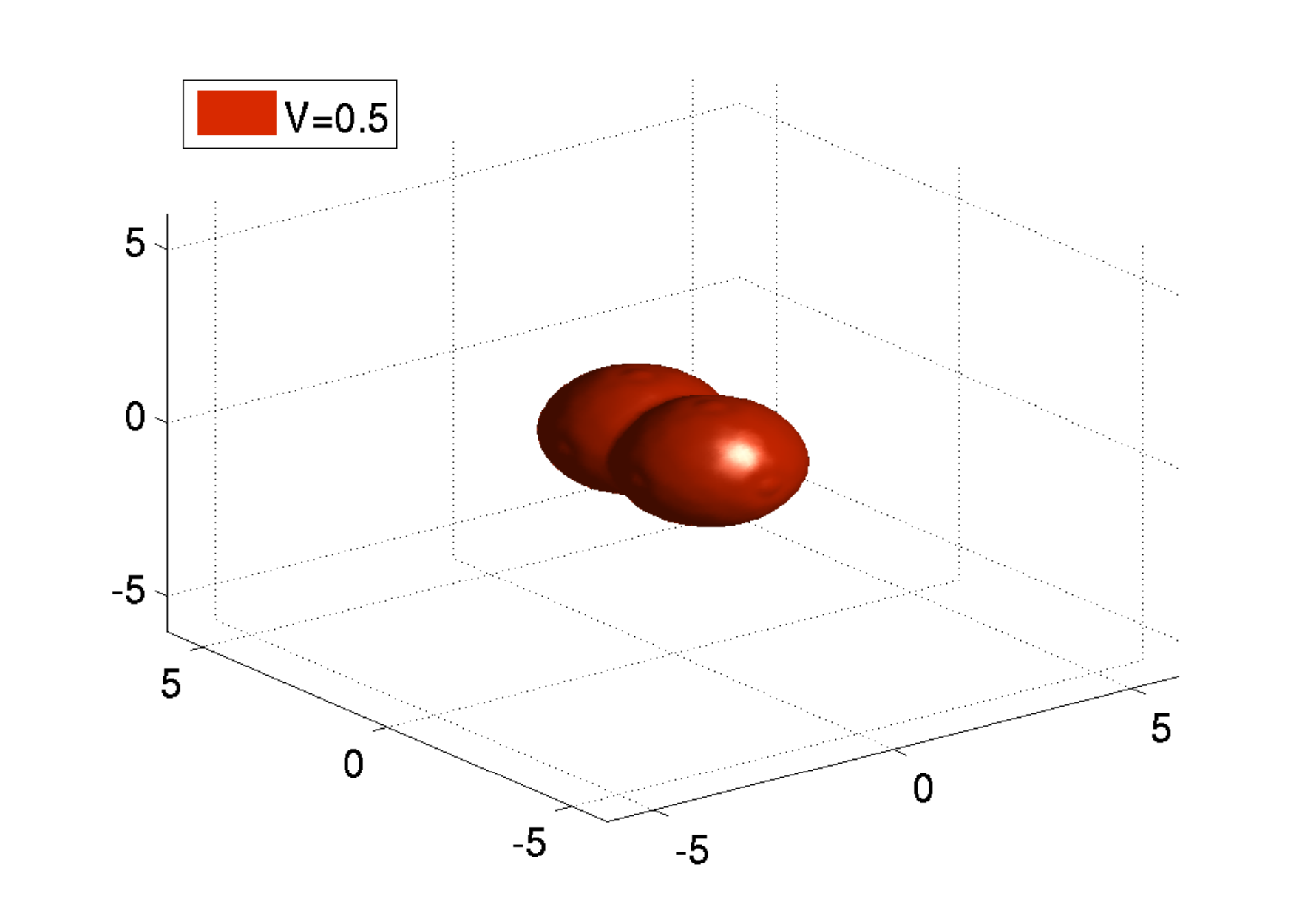,width=0.5\linewidth,clip=} &
\epsfig{file=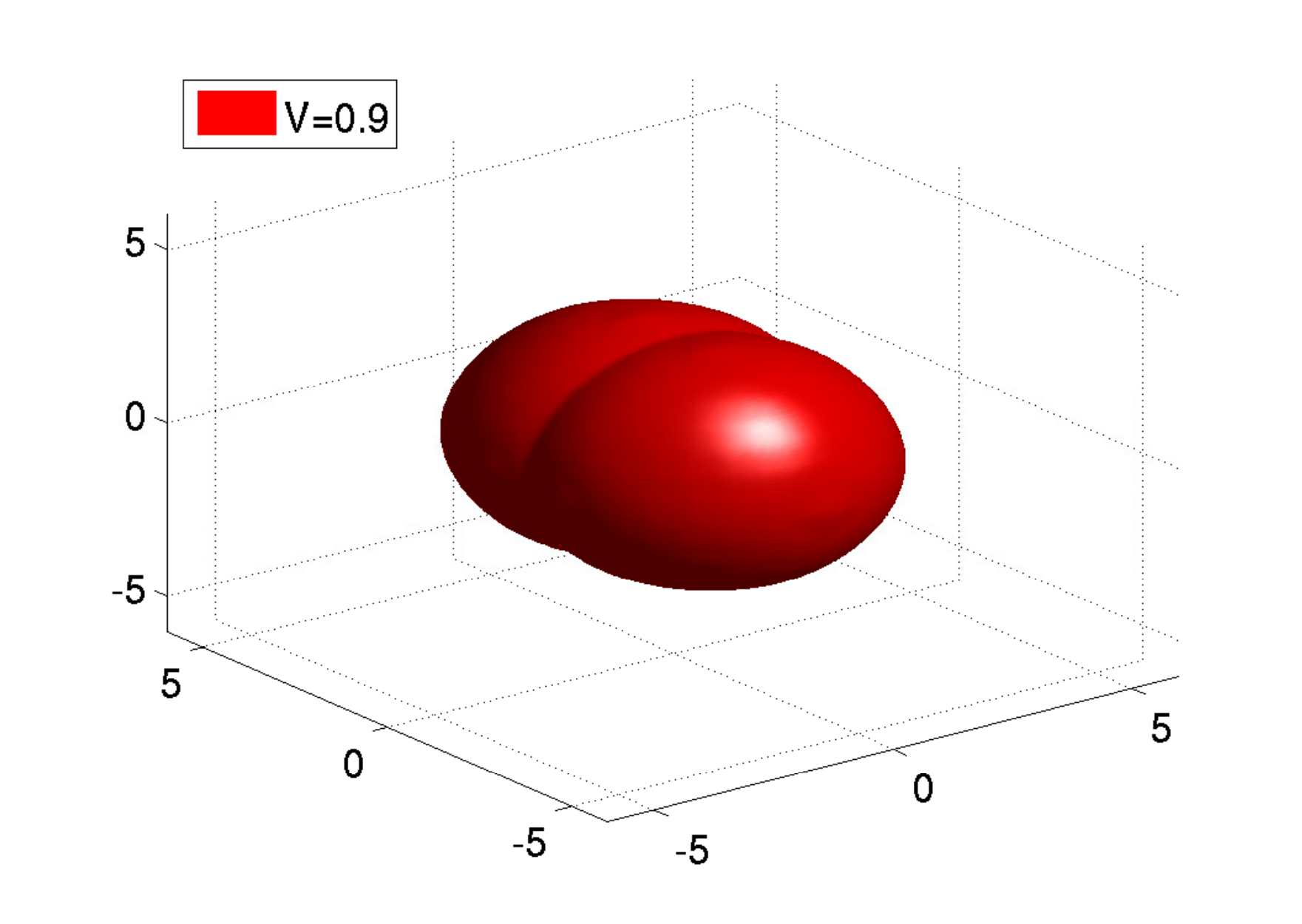,width=0.5\linewidth,clip=} 
\end{tabular}}
\caption{Test 7 (3D eikonal): different value function isosurfaces.}\label{fig:con3D}
\end{figure}

As in the two-dimensional study, we perform different tests by changing the domain and the target. For Test 6 we set $\Omega=]-1,\,1[^3$ and $\mathcal{T}=(0,0,0)$, while for Test 7, $\Omega=]-6,\,6[$ and $\mathcal{T}$ is the union of two unit spheres centered at $(-1,0,0)$ and $(1,0,0)$. In both cases, the set of controls is discretized into $16\times 8$ points. Reachable sets for Test 7 are shown in Figure \ref{fig:con3D}, and CPU times for both tests can be found in Tables \ref{tab:3D1} and \ref{tab:3D2}. We observe similar results as in the 2D tests, with up to $10\times$ acceleration for a simple target, and $4\times$ with more complicated targets. Note that in the second case, the speedup is similar to the natural performance that would be achieved by a PI algorithm. This is due to a weaker influence of the coarse VI iteration, which is sensitive to poor resolution of a complex target.

\begin{table}[ht]
\begin{center}
\ra{1.2}
\resizebox{.8\textwidth}{!}{
\begin{tabular}{cccccccccc}
\toprule
$\#$ nodes & $\Delta x$& &VI & &PI & & VI($2\Delta x$)& PI($\Delta x$) & API \\
\cmidrule{4-4}\cmidrule{6-6}\cmidrule{1-2}\cmidrule{8-10}
$41^3$       & 0.05        & & 4.83E2 (44)&   & 1.22E3 (10)  & &  4.61 (5)   & 1.19E2 (3)& 1.23E2 (8)\\
$81^3$      & 0.025      & & 7.67E3 (84)&   & 1.47E3 (13) & & 2.43E1 (12)  & 3.88E2 (3)& 6.31E2 (15)\\
\bottomrule
\end{tabular}}
\vskip 2mm
\caption{Test 6 (3D Eikonal): CPU time (iterations) for different algorithms with $a_1=16$ controls, $a_2=8$ controls.}\label{tab:3D1}
\end{center}
\end{table}

\begin{table}[ht]
\begin{center}
\ra{1.2}
\resizebox{.8\textwidth}{!}{
\begin{tabular}{cccccccccc}
\toprule
$\#$ nodes & $\Delta x$& &VI & &PI & & VI($2\Delta x$)& PI($\Delta x$) & API \\
\cmidrule{4-4}\cmidrule{6-6}\cmidrule{1-2}\cmidrule{8-10}
$61^3$       & 0.2        & & 2.67E2 (25)&   & 1.22E2 (9)  & &  1.44 (11)   & 6.80E1 (3)& 6.94E1 (14)\\
$121^3$      & 0.1      & & 4.52E3 (52)&   & 1.28E3 (11) & & 25.15E1 (12)  & 9.96E2 (3)& 1.01E3 (15)\\
\bottomrule
\end{tabular}}
\vskip 2mm
\caption{Test 7 (3D Eikonal): CPU time (iterations) for different algorithms with $a_1=16$ controls, $a_2=8$ controls.}\label{tab:3D2}
\end{center}
\end{table}

\subsection{Test 8: A Minimum time Problem in 4D}
\begin{figure}[ht]
\centering
\resizebox{0.7\textwidth}{!}{
\begin{tabular}{cc}
\epsfig{file=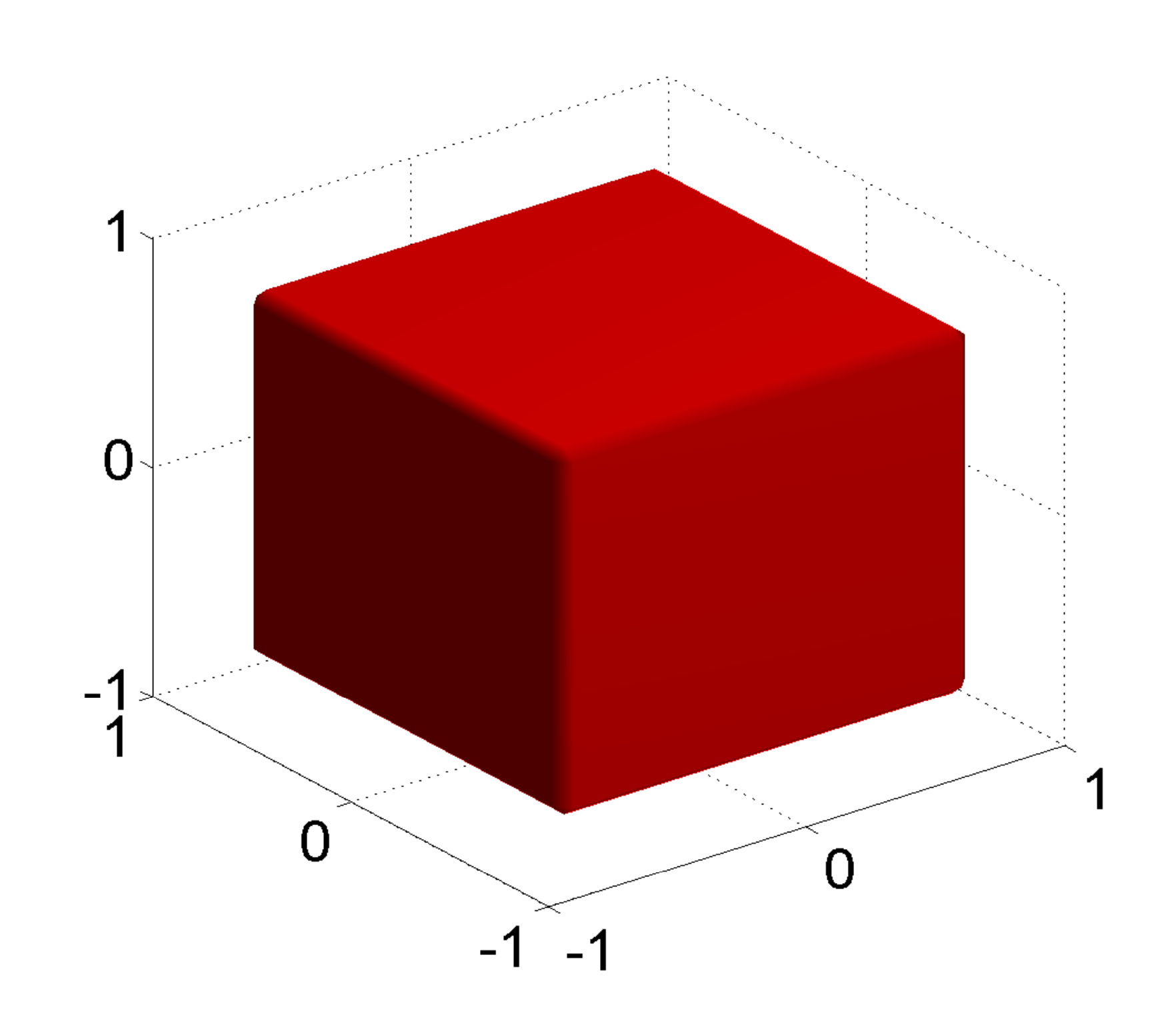,width=0.5\linewidth,clip=} &
\epsfig{file=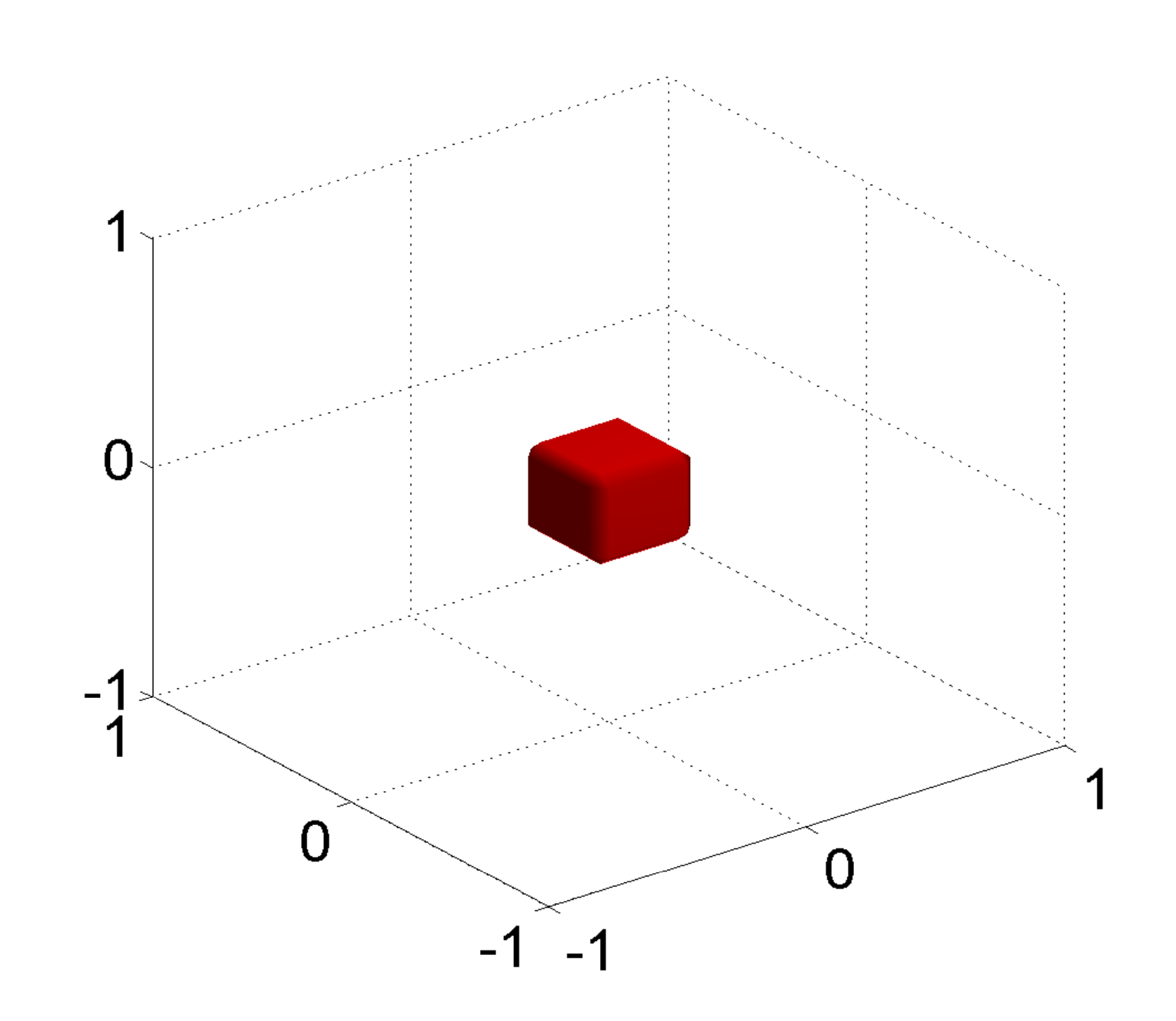,width=0.5\linewidth,clip=}
\end{tabular}}
\caption{Test 8 (4D minimum time): different value function isosurfaces with $x_4=0$.}
\label{fig:con4Dnew}
\end{figure}

We conclude our series of tests in minimum time problems by considering a four-dimensional problem with a relatively reduced control space. In the previous examples we have studied the performance of our scheme in cases where the set of discrete controls was fairly large, while in several applications, it is also often the case that the set of admissible discrete controls is limited and attention is directed towards the dimensionality of the state space. The following problem tries to mimic such setting. System dynamics are given by
\begin{equation*}
f(x,y,z,w,(a_1,a_2,a_3,a_4))=\left(\begin{array}{c}
a_1\\
a_2\\
a_3\\
a_4\end{array}\right)\,,
\end{equation*}
the domain is $\Omega=]-1,\,1[^4$, the target is $\mathcal{T}=\partial\Omega$, $\Delta t=0.8\Delta x$ and $A$ is the set of 8 directions pointing to the facets of the four-dimensional hypercube. Figure \ref{fig:con4Dnew} shows different reachable sets and CPU times are presented in Table \ref{tab:4D}. In the finest mesh a speedup of $8\times$ is observed, which is consistent with the previous results on simple targets. Thus, the performance of the presented algorithm is not sensitive neither to the number of discrete controls nor to the dimension of the state space, whereas it is affected by the complexity of the target.

\begin{table}[ht]
\begin{center}
\ra{1.2}
\resizebox{.8\textwidth}{!}{
\begin{tabular}{cccccccccc}
\toprule
$\#$ nodes & $\Delta x$& &VI & &PI & & VI($2\Delta x$)& PI($\Delta x$) & API \\
\cmidrule{4-4}\cmidrule{6-6}\cmidrule{1-2}\cmidrule{8-10}
$21^4$       & 0.1        & & 13.6 (15)&   & 16.2 (11)  & &  0.30 (4)  & 2.79 (2)& 3.09 (6)\\
$41^4$      & 5E-2     & & 4.79E2 (29)&   & 6.30E2 (21) & & 10.2 (12)  & 48.3 (2)& 58.5 (14)\\
\bottomrule
\end{tabular}}
\vskip 2mm
\caption{Test 8 (4D minimum time): CPU time (iterations) for different algorithms.}\label{tab:4D}
\end{center}
\end{table}

\subsection{Application to optimal control problem of PDE´s}
Having developed a comprehensive set of numerical tests concerning the solution of optimal control problems via static HJB equations, which assessed the performance of the proposed API algorithm, we present an application where the existence of accelerated solution techniques for high-dimensional problems is particularly relevant, namely, the optimal control of systems governed by partial differential equations.
From an abstract perspective, optimal control problems where the dynamics are given by evolutive partial differential equations correspond to systems where the state lies in an infinite-dimensional Hilbert space (see \cite{T10}). Nevertheless, in terms of practical applications, different discretization arguments can be used to deal with this fact, and (sub)optimal control synthesis can be achieved through finite-dimensional, large-scale approximations of the system. At this step, the resulting large-scale version will scale according to a finite element mesh parameter, and excepting for the linear-quadratic case and some closely related versions, it would be still computationally intractable for modern architectures (for instance, for a 100 elements discretization of a 1D PDE, the resulting optimal control would be characterized as the solution of a HJB equation in $\mathbb{R}^{100}$). Therefore, a standard remedy in optimal control and estimation is the application of model order reduction techniques, which, upon a large-scale version of the system, recover its most relevant dynamical features in a low-order approximation of prescribed size. In this context, surprisingly good control synthesis can be achieved with a reduced number of states (for complex nonlinear dynamics and control configurations an increased number of reduced states may be required). Previous attempts in this direction dates back to \cite{KVX04,KX05} and more recently to \cite{AF12,AF13bis}. We present an example where we embedded our accelerated algorithm inside the described framework. Note that, in this example, model reduction method is only applied in order to make the problem feasible for the Dynamic Programming approach. The acceleration is due to the  proposed API scheme.

\begin{figure}[ht]
\includegraphics[scale=0.4]{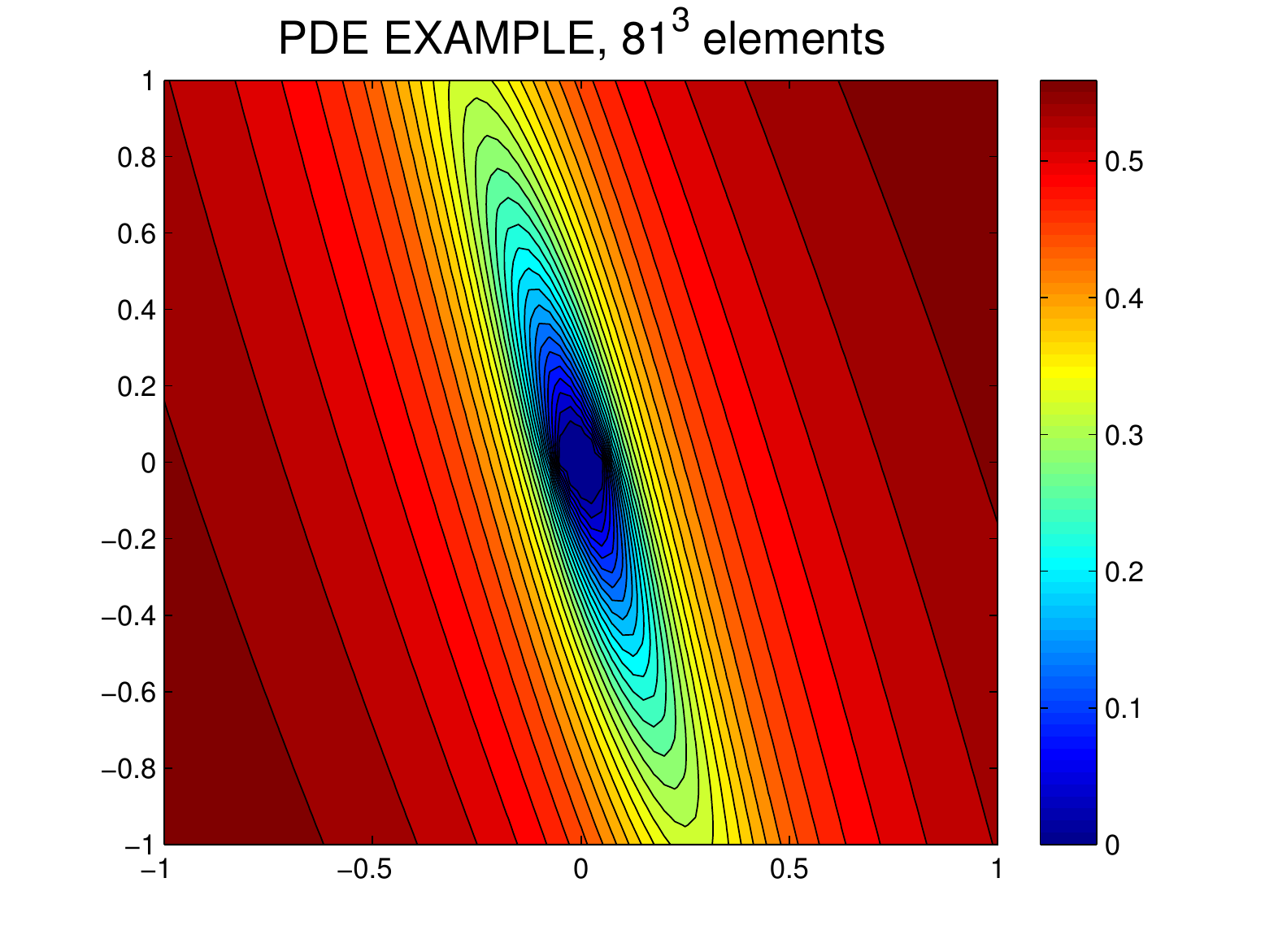}\hfill
\includegraphics[scale=0.4]{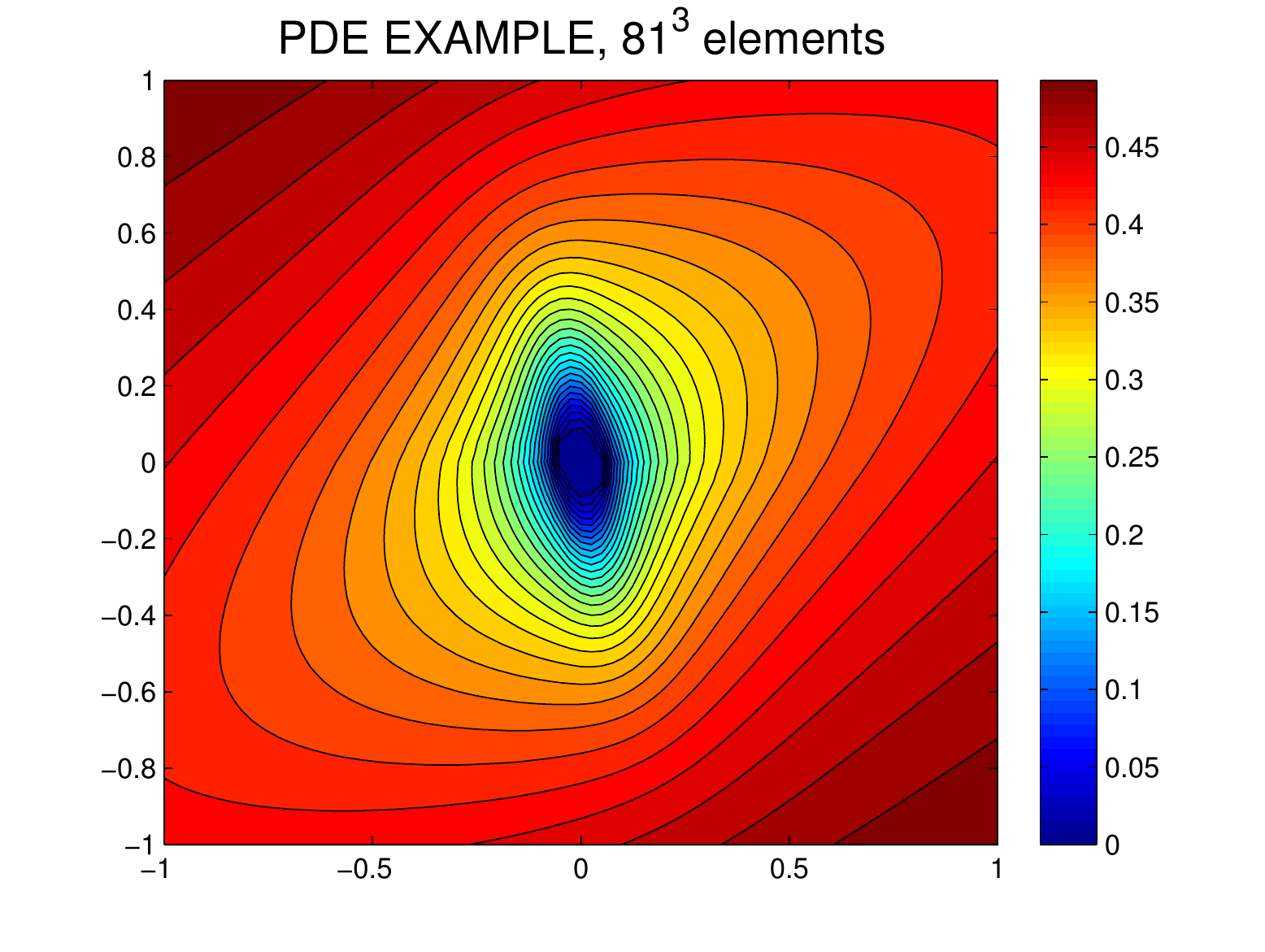}\hfill
\includegraphics[scale=0.4]{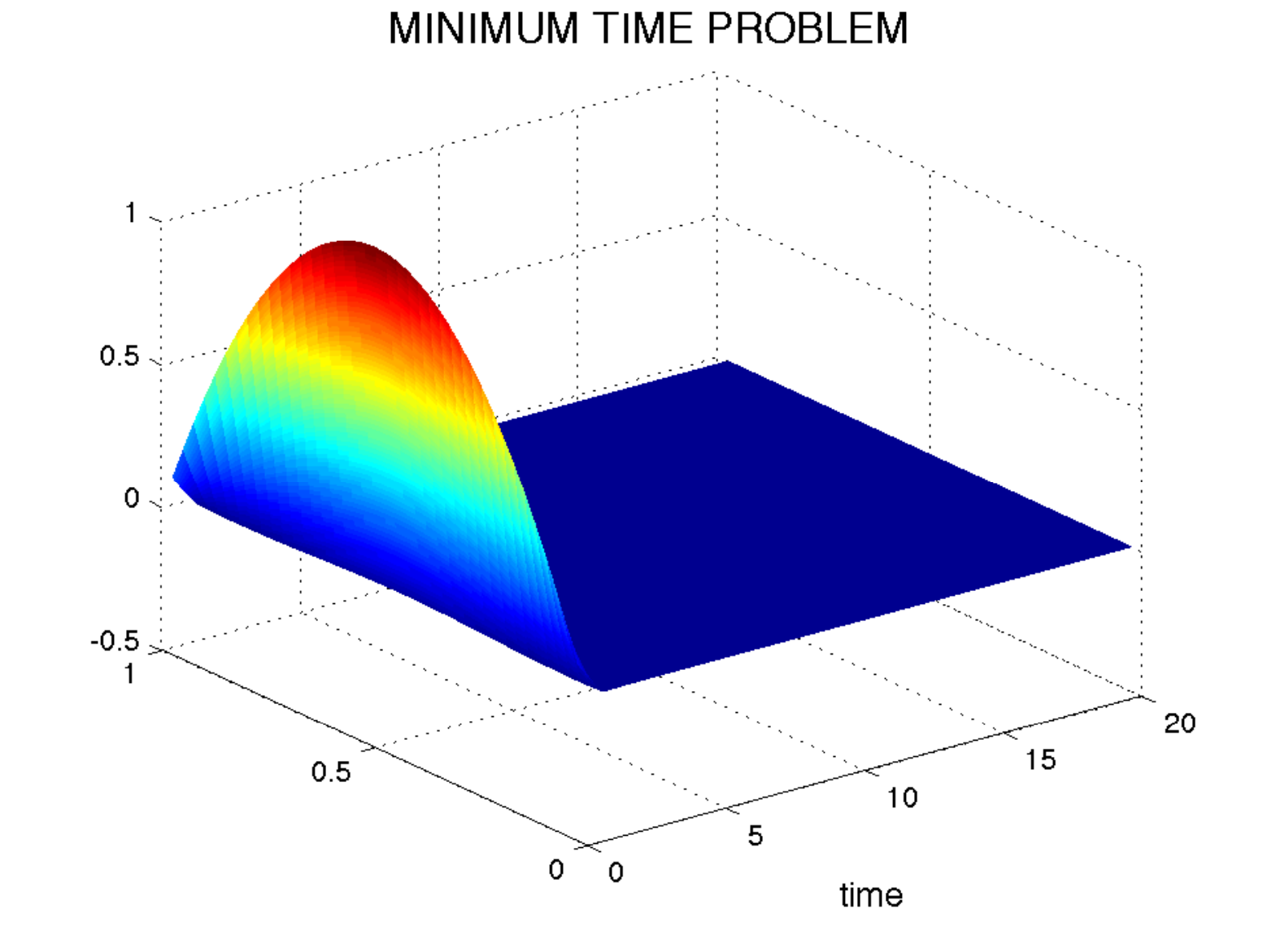}\hfill
\includegraphics[scale=0.4]{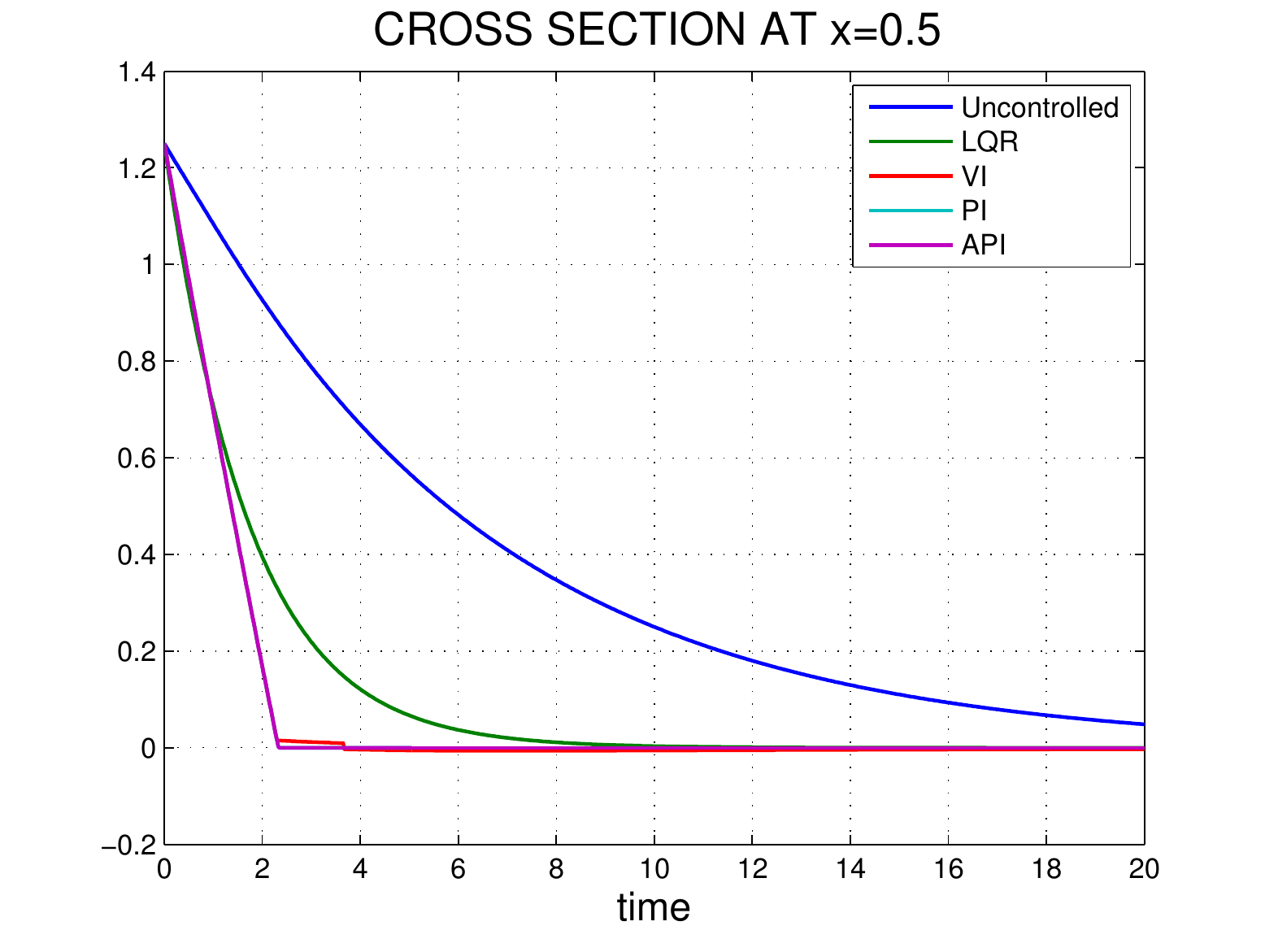}
\caption{Optimal control of the heat equation. Top left: contour plot of the value function at $x_3=0$. Top right: contour plot of the value function at $x_2=0$. Bottom left: controlled output via proposed procedure  of model reduction + minimum time HJB controller. Bottom right: cross section of the different outputs.}
\label{fig:heat}
\end{figure}

Let us consider a minimum time problem for the linear heat equation:
\begin{equation}\label{eq:heat}
\left\{\begin{array}{l}
y_t(x,t)=c y_{xx}(x,t)+y_0(x)\alpha(t)\,,\\
y(0,t)=y(1,t)=0\,,\\
y(x,0)=y_0(x)\,,
\end{array}\right.
\end{equation}
where $x\in[0,1],\,t\in [0, T]\,$, $c=1/80\,$ and  $\alpha(t):[0,T]\rightarrow\{-1,0,1\}$. After performing a finite difference discretization, we perform a Galerkin projection with basis function computed with a Proper Orthogonal Decomposition (POD) method, leading to a reduced order model (we refer to \cite{V11} for an introduction to this topic). In general,  model reduction techniques do have either a priori or a posteriori error estimates which allow to prescribe a certain number of reduced states yielding a desired level of accuracy. For this simple case, we consider the first 3 reduced states, which for a one-dimensional heat transfer process with one external source provides a reasonable description of the input-output behavior of the system. The system is reduced to:
\begin{equation}
\frac{d}{dt}\left[\begin{array}{c} x_1\\x_2\\x_3\end{array}\right]=\left[\begin{array}{ccc} -0.123&-0.008&-0.001\\-0.008&-1.148&-0.321\\-0.001&-0.321&-3.671\end{array}\right]\left[\begin{array}{c} x_1\\x_2\\x_3\end{array}\right]+\left[\begin{array}{c} -5.770\\-0.174\\-0.022\end{array}\right] \alpha(t).
\end{equation}
Once the reduced model has been obtained, we solve the minimum time problem with target $\mathcal{T}=(0,0,0)$. Figure \ref{fig:heat} shows contour plots of the value function in the reduced space and a comparison of the performance of the minimum time controller with respect to the uncontrolled solution and to a classical linear-quadratic controller is presented. CPU times are included in Table \ref{tablepde1}, where a speedup of $4\times$ can be observed, the acceleration would become more relevant as soon as more refined meshes and complex control configurations are considered.

\begin{table}[ht]
\begin{center}
\ra{1.2}
\resizebox{.8\textwidth}{!}{
\begin{tabular}{cccccccccc}
\toprule
$\#$ nodes & $\Delta x$& &VI & &PI & & VI($2\Delta x$)& PI($\Delta x$) & API \\
\cmidrule{4-4}\cmidrule{6-6}\cmidrule{1-2}\cmidrule{8-10}
$21^3$       & 0.1        & & 1.87 (76)&   & 0.91 (11)  & &  0.32 (27)  & 0.59 (8)& 0.98 (35)\\
$41^3$      & 5E-2      & & 27.8 (178)&  & 12.4 (15) & & 1.65 (76)  & 6.34 (10)& 7.99 (86)\\
$81^3$      & 2.5E-2      & & 6.13E2 (394)&  & 2.68E2 (15) & & 27.7 (178)  & 1.45E2 (9) & 1.72E2 (187)\\
\bottomrule
\end{tabular}}
\vskip 2mm
\caption{Minimum time control of the heat equation: CPU time (iterations) for different algorithms.}\label{tablepde1}
\end{center}
\end{table}


\section*{Concluding remarks and future directions}
In this work we have presented an accelerated algorithm for the solution of static HJB equations arising in different optimal control problems. The proposed method considers a pre-processing value iteration procedure over a coarse mesh with relaxed stopping criteria, which is used to generate a good initial guess for a policy iteration algorithm. This leads to accelerated numerical convergence with respect to the known approximation methods, with a speedup ranging in average from $4\times$ to $8\times$. We have assessed the performance of the new scheme via a extensive set of numerical tests focusing on infinite horizon and minimum time problems, providing numerical evidence of the reliability of the method in tests with increasing complexity. Positive aspects of the proposed scheme are its wide applicability spectrum (in general for static HJB), and its insensitivity with respect to the complexity of the discretized control set. Nonetheless, for some non trivial targets, special care is needed in order to ensure that the coarse pre-processing step will actually lead to an improved behavior of the policy iteration scheme.
Certainly, several directions of research remain open. The aim of this article was to present the numerical scheme and provide a numerical assessment of its potential. Future works should focus on tuning the algorithm in order to achieve an optimal performance; for instance, in order to make a fair comparison with the value iteration algorithm, the policy iteration step was also performed via a successive approximation scheme, while better results could be obtained by using a more efficient solver, including a larger amount of pre-processing work. Other possible improvements would relate to multigrid methods, high-order schemes and domain decomposition techniques. An area of application that remains unexplored is the case of differential games, where Hamilton-Jacobi-Isaacs equations need to be solved. Results presented in \cite{BMZ09} indicate that the extension is far from being trivial since a convergence framework is not easily guaranteed and the policy iteration scheme requires some modifications.


\end{document}